\documentclass[11pt]{article} 

\usepackage[utf8]{inputenc} 
\usepackage[section]{placeins} 

\usepackage{geometry} 
\usepackage{cite}
\usepackage{graphicx} 
\usepackage[parfill]{parskip} 

\usepackage{booktabs} 
\usepackage{array} 
\usepackage{paralist} 
\usepackage{verbatim} 
\usepackage{subfig} 
\usepackage{amssymb}
\usepackage{amsmath}
\usepackage{physics}
\usepackage{pifont}					
\usepackage{enumitem} 					

\usepackage[font=small,labelfont=bf,justification=raggedright,singlelinecheck=false]{caption}

\usepackage{fancyhdr} 
\usepackage{sectsty}
\usepackage{xcolor}
\definecolor{DarkGrey}{RGB}{105,105,105}
\definecolor{MidnightBlue}{RGB}{10,10,44}
\definecolor{RoyalBlue}{RGB}{25,41,88}
\definecolor{MyRoyalBlue}{RGB}{20,60,155}
\definecolor{GreyBlue}{RGB}{105,133,175}  
\definecolor{MyBlue}{RGB}{51,102,200}
\definecolor{LightBlue}{RGB}{0,153,255}
\definecolor{MyDarkBlue}{RGB}{0,20,100} 
\sectionfont{\color{MyBlue}}  
\subsectionfont{\color{MyRoyalBlue}}  
\subsubsectionfont{\color{MyDarkBlue}} 

\usepackage[nottoc,notlof,notlot]{tocbibind} 
\graphicspath{ {./figures/} }
\DeclareUnicodeCharacter{00B0}{$\,^{\circ}$}

\begin{document}

\title{Analytical exact gradients for single spin pulse sequence optimizations}
\author{Stella Slad and Burkhard Luy}
\date{} 

\maketitle

\begin{center}
{\it Institute for Biological Interfaces 4 – Magnetic Resonance, Karlsruhe Institute of Technology (KIT), Hermann-von-Helmholtz-Platz 1, 76344 Eggenstein-Leopoldshafen, Germany; 
Institute of Organic Chemistry, Karlsruhe Institute of Technology (KIT), Fritz-Haber-Weg 6, 76131 Karlsruhe, Germany}
\end{center}

\vspace{150px}





%

%

%
%

\begin{abstract}
The efficient computer optimization of magnetic resonance pulses and pulse sequences involves the calculation of a problem-adapted cost function as well as its gradients with respect to all controls applied. The gradients generally can be calculated as a finite difference approximation, as a GRAPE approximation, or as an exact function, e.g. by the use of the augmented matrix exponentiation, where the exact gradient should lead to best optimization convergence. However, calculation of exact gradients is computationally expensive and analytical exact solutions to the problem would be highly desirable. As the majority of todays pulse optimizations involve a single spin 1/2, which can be represented by simple rotation matrices in the Bloch space or by their corresponding Cayley-Klein/quaternion parameters, the derivations of analytical exact gradient functions appear to be feasible. Taking two optimization types, the optimization of point-to-point pulses using 3D-rotations and the optimization of universal rotation pulses using quaternions, analytical solutions for gradients with respect to controls have been derived. Controls in this case can be conventional $x$ and $y$ pulses, but also $z$-controls, as well as gradients with respect to amplitude and phase of a pulse shape. In addition, analytical solutions with respect to pseudo controls, involving holonomic constraints to maximum rf-amplitudes, maximum rf-power, or maximum rf-energy, are introduced. Using the hyperbolic tangent function, maximum values are imposed in a fully continuous and differentiable way. The obtained analytical gradients allow the calculation two orders of magnitude faster than the augmented matrix exponential approach. The exact gradients for different controls are finally compared in a number of optimizations involving broadband pulses for $^{15}$N, $^{13}$C, and $^{19}$F applications.
\end{abstract}

\newpage

\section{Introduction}
Computer optimization of pulses and pulse sequence building blocks has a long standing history in the NMR community. Starting out with composite pulses \cite{levittNMRPopulationInversion1979, levittCompositePulses1986, lurieNumericalDesignComposite1986, tyckoCompositePulsesPhase1985}, shaped pulses \cite{emsleyGaussianPulseCascades1990, emsleyOptimizationShapedSelective1992, ewingDevelopmentOptimizationShaped1990, garwoodSymmetricPulsesInduce1991, kupcePolychromaticSelectivePulses1993, kupceWidebandExcitationPolychromatic1994, smithImprovedBroadbandInversion2001, tannusAdiabaticPulses1997, zaxAmplitudemodulatedCompositePulses1988}, heteronuclear decoupling \cite{levittBroadbandHeteronuclearDecoupling1982, shakaEvaluationNewBroadband1983a} and Hartmann-Hahn-transfer building blocks \cite{kadkhodaieBroadbandHomonuclearCross1991, shakaIterativeSchemesBilinear1988}, the success of todays liquid state NMR spectroscopy is largely built on these elements. With the advent of optimal control based algorithms \cite{conollyOptimalControlSolutions1986, maoSelectiveInversionRadiofrequency1986, rosenfeldDesignAdiabaticSelective1996} and in particular with the GRAPE algorithm \cite{defouquieresSecondOrderGradient2011b, khanejaOptimalControlCoupled2005, skinnerApplicationOptimalControl2003b}, possibilities in pulse design increased dramatically in the past two decades and examples in pulse shape design \cite{ehniBEBEtrBUBIJcompensated2013b, ehniConcurrentJevolvingRefocusing2022a, ehniRobustINEPTRefocused2014, ehniSystematicApproachOptimizing2012, gershenzonLinearPhaseSlope2008, gershenzonOptimalControlDesign2007, hallerSORDORPulsesExpansion2022a, kobzarPatternPulsesDesign2005, koosBroadbandExcitationPulses2015, koosBroadbandRFAmplitude2017, martikyanApplicationOptimalControl2021b, sladBandselectiveUniversal902022a, spindlerShapedOptimalControl2012b, spindlerShapedOptimalControl2012b, spindlerPerspectivesShapedPulses2017b, odedraUseCompositeRefocusing2012, spindlerBroadbandInversionPELDOR2013, josephOptimalControlPulses2023, buchananSeedlessOntheflyPulse2024, heDigitalTwinParallel2024} show utmost performance close to the physical limits \cite{kobzarExploringLimitsBroadband2004, kobzarExploringLimitsBroadband2008a, kobzarExploringLimitsBroadband2012a, lapertExploringPhysicalLimits2012b}. However, even larger bandwidths and/or lower rf-energies will have to be explored as well as extremely complex optimizations like whole heternuclear decoupling periods \cite{nevesHeteronuclearDecouplingOptimal2009b, schillingNextGenerationHeteronuclearDecoupling2014c}, that easily lead to very long optimization times lasting weeks to months on highly parallelized supercomputers. It is therefore worth looking into the basic mathematics to look for analytical solutions of spin system treatment wherever possible to significantly speed up corresponding calculations. 

The majority of current pulse and pulse sequence optimizations are performed using a single spin in Liouville superoperator \cite{hogbenSpinachSoftwareLibrary2011b} or Bloch space \cite{bonnardReviewGeometricOptimal2012b, skinnerReducingDurationBroadband2004}. This description is perfectly suited for optimizing the transfer from an intial to a final state, as e.g. the case for so-called point-to-point pulses \cite{kobzarExploringLimitsBroadband2004, skinnerApplicationOptimalControl2003b,kobzarExploringLimitsBroadband2008a}. Other optimizations, like having a defined propagator  in universal rotation pulses \cite{khanejaOptimalControlCoupled2005, kobzarExploringLimitsBroadband2012a} as a target, are better described using either Cayley-Klein parameters \cite{cayleyIVSixthMemoir1997b} or the related quaternion formalism \cite{hamiltonQuaternionsNewSystem1843b, rodriguesLoisGeometriquesQui1840b}.  In each type of optimization, a Hamiltonian needs to be diagonalized, which can be done by numerical exponentiation using for example the Padé algorithm, or by analytical formulae, which are particularly straightforward in the case of a single spin in Bloch space, where this is resembled by well-studied three dimensional Cartesian rotation matrices. In addition, partial derivatives to all controls need to be calculated in gradient and hessian based optimization algorithms. This can be achieved by a linear approximation as in the original GRAPE formalism \cite{khanejaOptimalControlCoupled2005, skinnerApplicationOptimalControl2003b}, but exact gradients are more efficient for the convergence of optimizations \cite{defouquieresSecondOrderGradient2011b, kuprovOptimalControlSpin2023b}. Such exact gradients can be calculated using a general approach described by Ernst \cite{aizuParameterDifferentiationQuantum1963b, levantePulseSequenceOptimizationAnalytical1996b}, or by a recent development in algebra based on the matrix exponentiation of co-propagators \cite{goodwinAuxiliaryMatrixFormalism2015, goodwinModifiedNewtonRaphsonGRAPE2016}. For a single spin, however, analytical solutions with respect to conventional controls were recently introduced with the ESCALADE approach \cite{foroozandehOptimalControlSpins2021} and it was shown that very fast optimizations are possible based on a mixture of exact analytical gradients and numerical hessian calculations \cite{goodwinAcceleratedNewtonRaphsonGRAPE2023b}. 

The goal of computational optimization in NMR spectroscopy must be to be able to address more and more complex problems. In this publication we focus on the improvement of extremely complex optimizations that can be performed with single spin calculations. One goal is, for example, to achieve broadband pulses that can cover bandwidths 20 to 100 times the applied maximum rf-amplitude to address the full chemical shift bandwidths of e.g. $^{15}$N, $^{19}$F, $^{31}$P, $^{119}$Sn, or $^{195}$Pt. In such cases, computational demands rise non-linearly and while short pulse shapes can be optimized in a matter of seconds to minutes on a conventional PC, the desired demanding pulses may take months on high level supercomputers or are not at all feasible with current technology due to limitations in memory. This memory issue is particularly pressing in the case of hessian-based optimization algorithms. If, for example, a 10~ms pulse with 0.1~$\mu$s digitization is considered, 100.000 digits result in 200.000 controls (x and y component), which in turn lead to a 200.000$\times$200.000 hessian matrix with double precision numbers, which corresponds to 320~GB memory for the calculation of a single condition hessian in an optimization. As our aim is at this type of optimization complexity, we do not consider second derivative hessian matrices in the following and reduce our view to exact gradients only, which can be calculated much more memory efficient and are readily combined with e.g. the LBFGS optimization algorithm. In this context, we here derive analytical derivatives in the Bloch picture based on rotational matrices, as well as for four-dimensional quaternions for propagator-based optimizations in a tabular form. It is thereby closely related to the ESCALADE approach, for which a very nice abstract way of gradient and hessian calculations based on the Rodrigues expansion has been developed \cite{foroozandehOptimalControlSpins2021}. In the present approach, however, we use a much simpler method in which the quaternions and the Rodrigues formula for rotations are directly used to derive gradients simply from individual matrix components. We therefore can extent the ESCALADE approach to more individualized optimization concepts: next to the conventional case with $x$ and $y$-controls, we also consider several special cases, like the formulation in Cartesian coordinates for $x$,$y$ pulses including $z$-controls, or the equivalent description using polar amplitude, phase and $z$-controls for pulse shapes. Also direct derivatives for holonomic constraints using $\tanh$ based pseudo-parameters for amplitude, power, and energy-limited pulse optimizations are derived.


\section{Theory}
\subsection{GRAPE algorithm formulation for a single spin $\frac{1}{2}$}
\subsubsection{Optimization of point-to-point pulses}

A shaped pulse can be seen as a sequence of $N$ short pulses of length $\Delta t$ with piecewise constant rf-amplitudes, where the $j^{\mathrm{th}}$ pulse normally consists of two controls $\omega_{\mathrm{x}}(j)$ and $\omega_{\mathrm{y}}(j)$ that represent the $x$ and $y$-components of the shaped pulse. While the offset frequency $\omega_{\mathrm{off}}$ at a certain position in the spectrum relative to the irradiation frequency constitutes the free evolution or drift Hamiltonian

\begin{equation}
{\cal{H}}_0 =  \omega_{\mathrm{off}} \ I_z = 2 \pi \ \nu_{\mathrm{off}} \  I_z,
\end{equation}

the pulses constitute the control Hamiltonian at the $j^{\mathrm{th}}$ pulse

\begin{equation}
{\cal{H}}_{\mathrm{1}}(j) =  \omega_{\mathrm{x}}(j) \ I_x + \omega_{\mathrm{y}}(j) \ I_y = 2 \pi \ \nu_{\mathrm{rf}}(j)  \{ \cos \alpha(j) \  I_x + \sin \alpha(j) \  I_y \}.
\end{equation}

For a given sequence of $N$ pulses with duration $\Delta t$ at a specific offset $\nu_{\mathrm{off}}$ the propagator is given by 

\begin{equation}
U_j = \exp \large\{ - i \{ {\cal{H}}_0 +  {\cal{H}}_{\mathrm{1}}(j) \} \Delta t \large \}
\end{equation}

and the propagation of an initial spin density operator $\rho_0$ can be written as 

\begin{equation}
\rho_N = U_N \cdots U_j \cdots U_1 \ \rho_0 \ U^{\dagger}_1 \cdots U^{\dagger}_j \cdots U^{\dagger}_N . 
\end{equation}


The goal of a point-to-point optimization is to find values of the controls that minimize the differences between the desired final state of the spin ($\mathrm{\lambda}_{\mathrm{F}}$) and the obtained final density operator with the current control amplitudes ($\mathrm{\rho}_N$), which is identical to maximizing their overlap according to \cite{khanejaOptimalControlCoupled2005}

\begin{equation}
\Phi_{\mathrm{PP}}= Re \braket{\lambda_{\mathrm{F}}}{\rho_{{N}}}
\end{equation}


In order to maximize the cost function $\Phi_{\mathrm{PP}}$, we have to minimize its gradient with respect to every control at every timestep, which for the $j^{th}$ time point  is resulting in

\begin{equation}
\label{eq:gradient_pp}
\Gamma_{\mathrm{PP}}(j)=\begin{pmatrix}\dfrac{\partial\Phi_{\mathrm{PP}}}{\partial\omega_x(j)} \\[4mm] \dfrac{\partial\Phi_{\mathrm{PP}}}{\partial\omega_y(j)}
\end{pmatrix}
= \begin{pmatrix}\braket{\lambda_j}{ \dfrac{\partial U_j}{\partial\omega_x(j)} \rho_{j-1} U^{\dagger}_j} 
+ \braket{\lambda_j}{  U_j \rho_{j-1} \dfrac{\partial U^{\dagger}_j}{\partial\omega_x(j)}}\\[4mm]  \braket{\lambda_j}{ \dfrac{\partial U_j}{\partial\omega_y(j)} \rho_{j-1} U^{\dagger}_j} 
+ \braket{\lambda_j}{  U_j \rho_{j-1} \dfrac{\partial U^{\dagger}_j}{\partial\omega_y(j)}} 
\end{pmatrix}
\approx 
\begin{pmatrix}\braket{\lambda_j}{ - i \Delta t [{\cal{H}}_{1_x}(j), \rho_j]} \\[4mm]  \braket{\lambda_j}{ - i \Delta t [{\cal{H}}_{1_y}(j), \rho_j]} 
\end{pmatrix}
\end{equation}

with $\lambda_j = U^{\dagger}_{j+1} \cdots U^{\dagger}_{N} \ \lambda_{\mathrm{F}} U_N \cdots U_{j+1}$, $\rho_j =  U_j \cdots U_1 \ \rho_0 \ U^{\dagger}_1 \cdots U^{\dagger}_j$ and ${\cal{H}}_{1_x}(j)$ and ${\cal{H}}_{1_y}(j)$ being the $x$ and $y$ components of the control Hamiltonian.   
Thus, for each iteration $i$ of the optimization, the controls $\omega^{\mathrm{}}_{{k}}(j)$ , $k\in \{x,y\}$,
are guaranteed to increase the cost function $\Phi_{\mathrm{PP}}$ for an infinitesimal $\epsilon$ according to 
\begin{equation} \label{eq:update_controls}
\omega_{{k}}^{(i+1)}(j) \rightarrow \omega_{{k}}^{(i)}(j)+\epsilon \frac{\partial \Phi_{\mathrm{PP}}}{\partial \omega_{{k}}(j)},
\end{equation}
until convergence to a (local) optimum is reached. 
Please note that the partial derivatives may also be calculated using the rotation angles 
\begin{equation}
\theta_x(j) = \omega_x(j) \, \Delta t \, \, ; \, \, \theta_y(j) = \omega_y(j) \, \Delta t.
\end{equation}
In this case, the update can be rewritten into
\begin{equation} \label{eq:update_controls_theta}
\theta_{{k}}^{(i+1)}(j) \rightarrow \theta_{{k}}^{(i)}(j)+\epsilon \frac{\partial \Phi_{\mathrm{PP}}}{\partial \theta_{{k}}(j)},
\end{equation}
using the relation $d\omega / d\theta = 1/\Delta t $.  
Instead of the GRAPE-approximation given at the end of Eq. \ref{eq:gradient_pp}, it is advantageous to use the computationally more costly exact gradients for the pulse sequence update \cite{defouquieresSecondOrderGradient2011b}. The critical point in these exact gradients is the calculation of $\partial U_j / \partial\omega_x(j)$ or $\partial U_j / \partial\theta_x(j)$, respectively, and its complex conjugate, for which several solutions have been proposed \cite{goodwinAuxiliaryMatrixFormalism2015, levantePulseSequenceOptimizationAnalytical1996b}, but elegant analytical solutions might be of particular interest for specific problems. For a single spin we will derive such very compact analytical forms for various optimization types in the following.


\noindent
In order to design broadband PP pulses with $B_1$ inhomogeneity compensation, it is further necessary to calculate the cost function and the gradient for $n_{\mathrm{off}}$ offsets linearly distributed over the desired bandwdith $\Delta \nu$ and $n_{\mathrm{rf}}$ different rf-amplitudes $\nu_{\mathrm{rf}}$ covering the desired $B_1$-compensated range. The global quality factor and the global gradient can be calculated as averages according to 
\begin{equation} \overline{\Phi_{\mathrm{PP}}}=\frac{1}{n_{\mathrm{off}}n_{\mathrm{rf}}}\displaystyle\sum_{i=1}^{n_{\mathrm{off}}}\sum_{l=1}^{n_{\mathrm{rf}}} \Phi_{\mathrm{PP}}(\nu^i_{\mathrm{off}},\nu^l_{\mathrm{rf}}) \end{equation}

and

\begin{equation} \overline{\Gamma_{\mathrm{PP}}(j)}=\frac{1}{n_{\mathrm{off}}n_{\mathrm{rf}}}\sum_{i=1}^{n_{\mathrm{off}}}\sum_{l=1}^{n_{\mathrm{rf}}} \Gamma_{\mathrm{PP}}(j,\nu^i_{\mathrm{off}},\nu^l_{\mathrm{rf}}) .\end{equation}

\subsubsection{Optimization of universal rotation pulses}

If not the transformation of a given state to a final state, but rather the rotation in Hilbert space itself is the target, the propagator itself can be optimized \cite{defouquieresSecondOrderGradient2011b, khanejaOptimalControlCoupled2005, kobzarExploringLimitsBroadband2012a, luyConstructionUniversalRotations2005, skinnerNewStrategiesDesigning2012}. Considering the total propagator $U(T)$ at time point $T = N \Delta t$ given by

\begin{equation}
U(T) = U_N \cdots U_j \cdots U_1,
\end{equation}

a cost function with respect to the desired propagator $U_{\mathrm{F}}$ can be formulated as 

\begin{equation}
\Phi_{\mathrm{UR}} = Re\braket{U_{\mathrm{F}}}{U(T)},
\end{equation}

and the corresponding gradients to the $j^{th}$ time step for $x$ and $y$ components of the pulse train s are approximated by \cite{khanejaOptimalControlCoupled2005}
 
\begin{equation}
\label{eq:gradient_ur}
\Gamma_{\mathrm{UR}}(j)=\begin{pmatrix}\dfrac{\partial\Phi_{\mathrm{UR}}}{\partial\theta_x(j)} \\[4mm] \dfrac{\partial\Phi_{\mathrm{UR}}}{\partial\theta_y(j)}
\end{pmatrix}
= \begin{pmatrix}Re \braket{P_j}{\dfrac{\partial U_j}{\partial\theta_x(j)} X_{j-1}} \\[4mm]  
Re \braket{P_j}{\dfrac{\partial U_j}{\partial\theta_y(j)} X_{j-1}} 
\end{pmatrix}
\approx 
\begin{pmatrix}Re \braket{P_j}{ - i \Delta t [{\cal{H}}_{1_x}(j), X_j]} \\[4mm]  Re \braket{P_j}{ - i \Delta t [{\cal{H}}_{1_y}(j), X_j]} 
\end{pmatrix}
\end{equation}

with  $P_j = U^{\dagger}_{j+1} \cdots U^{\dagger}_{N} \ U_{\mathrm{F}}$, $X_j =  U_j \cdots U_1 $ and previously used notations for the control Hamiltonian components ${\cal{H}}_{1_x}(j)$ and ${\cal{H}}_{1_y}(j)$.

\subsection{Analytical exact gradients $\mathbf{\Gamma_{\mathrm{PP}}}$ using $\mathbf{x}$, $\mathbf{y}$, and optional $\mathbf{z}$ controls}

The usual equations for evolution of a spin density matrix have been used in the previous section. However, it might be advantageous to use the Liouville superoperator formalism for a single spin $\frac{1}{2}$ in Cartesian coordinate representation. The spin density operator can then be represented by a four vector $\rho = (\rho_{ \bf{1}}, \rho_x, \rho_y, \rho_z)^T$, where the contribution of the identity matrix may be neglected. If furthermore relaxation is neglected, the spin density at time point $j$ is fully represented by the vector 

\begin{equation}
\rho_{\mathrm{j}}=
\begin{pmatrix}
\rho_{\mathrm{x}} \\ \rho_{\mathrm{y}} \\ \rho_{\mathrm{z}}
\end{pmatrix}
\end{equation}

and propagation is achieved by the Liouville superoperators $R_j$ in the reduced Cartesian representation

\begin{equation} \label{eq:propagation}
\rho_{\mathrm{j}}=R_{\mathrm{j}}...R_{\mathrm{1}}\rho_0
\end{equation}

and the co-state is represented accordingly by

\begin{equation}
\lambda_j=R_{\mathrm{j}}^{-1}...R_{\mathrm{N}}^{-1}\lambda_{\mathrm{F}}.
\end{equation}

Please note, that the spin density as well as the co-state vector only have to be multiplied single-sided by the corresponding Liouville superoperators. Since we are in the Cartesian component basis set, we furthermore see that the superoperator is represented by a simple rotation matrix, which can be written as 

\begin{flalign}
\label{eq16}
R_{\mathrm{j}}= 
 \begin{pmatrix}
 \cos(\theta)+n_x^2(1-\cos(\theta)) & 
-n_z\sin(\theta)+n_xn_y(1-\cos(\theta)) & 
n_y\sin(\theta)+n_xn_z(1-\cos(\theta)) \\
n_z\sin(\theta)+n_xn_y(1-\cos(\theta)) &
\cos(\theta)+n_y^2(1-\cos(\theta))  &
-n_x\sin(\theta)+n_yn_z(1-\cos(\theta)) \\
-n_y\sin(\theta)+n_xn_z(1-\cos(\theta)) &
n_x\sin(\theta)+n_yn_z(1-\cos(\theta)) &
\cos(\theta)+n_z^2(1-\cos(\theta)),
 \end{pmatrix} 
 \end{flalign}
 with the overall rotation angle $\theta$ and the normalized components
$n_x$, $n_y$ and $n_z$  of the rotation axis of timestep $j$. It can also be noted more compact using the 
 Rodrigues formula using its individual elements $R_{hk}$:
 \begin{equation}
 \label{eq17}
 R_{hk}=\begin{cases}
     \cos^{2}\bigg(\dfrac{\theta}{2}\bigg)+(2n_h^2-1)\sin^{2}\bigg(\dfrac{\theta}{2}\bigg) & \text{für}~h=k \\
     2n_hn_k\sin^{2}\bigg(\dfrac{\theta}{2}\bigg)-\epsilon_{hkl}~n_l\sin(\theta) & \text{für}~h\neq k,
   \end{cases}
 \end{equation}
where $\epsilon_{hkl}$ is the Levi-Civita symbol. 
 
In each step $j$, $\theta$, $n_x$, $n_y$ and $n_z$ can be calculated from the control values $\omega_x$, $\omega_y$ and $\omega_z$ as follows:
\begin{flalign}
\label{eqtheta}
\theta &=\Delta t\sqrt{\omega_x^2+\omega_y^2+(\omega_z+\omega_{\mathrm{off}})^2} 
= \sqrt{\theta_x^2+\theta_y^2+\theta_z^2}  \\
n_x &=\dfrac{\omega_x}{\sqrt{\omega_x^2+\omega_y^2+(\omega_z+\omega_{\mathrm{off}})^2}}=\dfrac{\theta_x}{\sqrt{\theta_x^2+\theta_y^2+\theta_z^2}} = \dfrac{\theta_x}{\theta}
\\
n_y &= \dfrac{\omega_y}{\sqrt{\omega_x^2+\omega_y^2+(\omega_z+\omega_{\mathrm{off}})^2}}
=\dfrac{\theta_y}{\sqrt{\theta_x^2+\theta_y^2+\theta_z^2}} = \dfrac{\theta_y}{\theta}
\\
\label{eqthetaend}
n_z &= \dfrac{\omega_z+\omega_{\mathrm{off}}}{\sqrt{\omega_x^2+\omega_y^2+(\omega_z+\omega_{\mathrm{off}})^2}} 
=\dfrac{\theta_z+\theta_{\mathrm{off}}}{\sqrt{\theta_x^2+\theta_y^2+\theta_z^2}} = \dfrac{\theta_z}{\theta}.
\end{flalign}

with the effective flip angles around the Cartesian axes $\theta_x = \omega_x \Delta t$, $\theta_y = \omega_y \Delta t$, $\theta_z = (\omega_z + \omega_{\mathrm{off}}) \Delta t$, where the rotation around the $z$-axis now consists of the offset term introduced above and the introduction of a potential $z$-control that cannot be applied directly on a spectrometer, but that can be used to allow direct implementation of pulse sequence bound offset changes in optimizations \cite{stockmannVivoB0Field2018, vindingLocalSARGlobal2017, vindingOptimalControlGradient2021}. With these equations in hand, it is straightforward to reformulate the cost function for point-to-point pulses in Cartesian Liouvielle representation as

\begin{equation}
\Phi_{\mathrm{PP}}=\braket{\lambda_{\mathrm{F}}}{\rho_{\mathrm{N}}}=\braket{\lambda_{\mathrm{F}}}{R_{\mathrm{N}}...R_{\mathrm{1}}\rho_0} =\Big\langle\underbrace{R_{\mathrm{j+1}}^{T}...R_{\mathrm{n}}^{T}\lambda_{\mathrm{n}}}_{\textstyle\lambda_{\mathrm{j}}}\Big |
\underbrace{R_{\mathrm{j}}...R_{\mathrm{1}}\rho_0}_{\textstyle\rho_{\mathrm{j}}}\Big\rangle.
\end{equation}

Accordingly, the gradients for $j^{th}$ time point for the single spin is

\begin{equation}
\label{eq:gradient_ppR}
\Gamma_{\mathrm{PP}}(j)=\begin{pmatrix}\dfrac{\partial\Phi_{\mathrm{PP}}}{\partial\theta_x(j)} \\[4mm] \dfrac{\partial\Phi_{\mathrm{PP}}}{\partial\theta_y(j)} \\[4mm]
\dfrac{\partial\Phi_{\mathrm{PP}}}{\partial\theta_z(j)}
\end{pmatrix}
= \begin{pmatrix}\braket{\lambda_j}{ \dfrac{\partial R_j}{\partial\theta_x(j)} \rho_{j-1}}\\[4mm]  \braket{\lambda_j}{ \dfrac{\partial R_j}{\partial\theta_y(j)} \rho_{j-1}}\\[4mm] 
\braket{\lambda_j}{ \dfrac{\partial R_j}{\partial\theta_z(j)} \rho_{j-1}}
\end{pmatrix},
\end{equation}

which leaves an analytical derivation of the derivative of the rotation matrix with respect to the Cartesian controls to obtain an overall exact gradient. With the formulae given in this section, the derivatives of the individual matrix components can be straightforwardly achieved and corresponding results are shown in Table \ref{tab:pp_dx} for the $x$ and$y$ derivatives and in Table \ref{tab:pp_dz} for the $z$-control.




\subsection{Analytical exact gradients $\mathbf{\Gamma_{\mathrm{PP}}}$ using amplitude, phase, and optional $z$ controls}

The normalized components of the rotation axes $n_x$, $n_y$ and $n_z$ can also be written in polar coordinates, resulting in the pulse phase $\alpha$ and rf-amplitude $\omega_{\mathrm{rf}}$ as the input parameters for pulse shapes. The corresponding notation is
\begin{equation}
\label{eq24}
n_x =\dfrac{\cos(\alpha)\theta_{xy}}{\theta}=\dfrac{\cos(\alpha)\omega_{\mathrm{rf}}\Delta t}{\sqrt{(\Delta t^2\omega_{\mathrm{rf}})^2+\theta_z^2}},
~n_y =\dfrac{\sin(\alpha)\theta_{xy}}{\theta}=\dfrac{\sin(\alpha)\omega_{\mathrm{rf}}\Delta t}{\sqrt{(\Delta t^2\omega_{\mathrm{rf}})^2+\theta_z^2}},
~n_z =\dfrac{\theta_z}{\theta}=\dfrac{\theta_z}{\sqrt{(\Delta t^2\omega_{\mathrm{rf}})^2+\theta_z^2}}
\end{equation}
with 
\begin{flalign}
\theta_{xy} &=\omega_{\mathrm{rf}}\Delta t \notag\\
\theta_z &=(\omega_z+\omega_{\mathrm{off}})\Delta t  \notag \\
\theta &=\sqrt{\theta_{xy}^2 + \theta_z^2}
\label{eq25}
\end{flalign}

The derivatives for the rotational matrix components with respect to $\alpha$, $\omega_{\mathrm{rf}}$ and $\omega_z$ can then be derived the same way as in the previous section.  The gradients for the $j^{th}$ time point for the single spin are

\begin{equation}
\label{eq:gradient_ppR}
\Gamma_{\mathrm{PP}}(j)=\begin{pmatrix}\dfrac{\partial\Phi_{\mathrm{PP}}}{\partial \alpha(j)} \\[4mm] \dfrac{\partial\Phi_{\mathrm{PP}}}{\partial\theta_{\mathrm{xy}}(j)} \\[4mm]
\dfrac{\partial\Phi_{\mathrm{PP}}}{\partial\theta_z(j)}
\end{pmatrix}
= \begin{pmatrix}\braket{\lambda_j}{ \dfrac{\partial R_j}{\partial\alpha(j)} \rho_{j-1}}\\[4mm]  \braket{\lambda_j}{ \dfrac{\partial R_j}{\partial\theta_{\mathrm{xy}}(j)} \rho_{j-1}}\\[4mm] 
\braket{\lambda_j}{ \dfrac{\partial R_j}{\partial\theta_z(j)} \rho_{j-1}}
\end{pmatrix},
\end{equation}

which again leaves an analytical derivation of the derivative of the rotation matrix with respect to the Cartesian controls to obtain an overall exact gradient. The calculation can be done by hand or a symbolic mathematics program and the result for the rotation matrix components are listed in Tables~\ref{tab:pp_awz1} and~\ref{tab:pp_awz2}, respectively. The amplitude and phase representation is particularly useful in cases of restricted rf-amplitudes or even constant rf-amplitudes \cite{skinnerOptimalControlDesign2006}, where in the latter case only derivatives to the phase $\alpha$ need to be considered for a minimum set of parameters.

\FloatBarrier
\subsection{Analytical exact gradients $\mathbf{\Gamma_{\mathrm{UR}}}$}


We have seen that propagation can be expressed in terms of simple rotations in the case of a single spin $\frac{1}{2}$. If effective rotations themselves need to be optimized, it is best to express them with minmum storage and computation time. As such, Caley-Klein \cite{cayleyIVSixthMemoir1997b} or, equivalently, quternions \cite{emsleyGaussianPulseCascades1990, hamiltonQuaternionsNewSystem1843b, rodriguesLoisGeometriquesQui1840b, sladBandselectiveUniversal902022a} can be used to express rotations. Although the minimum set of numbers to represent a rotation is the vector spanned by $(\theta_x, \theta_y, \theta_z)$, it is better to use a four vector

\begin{equation}
\label{eq27}
Q_j =  
\begin{pmatrix}
A_j \\
B_j\\
C_j \\
D_j
\end{pmatrix}
\end{equation}

for the rotation at time step $j$, where the four components are defined by 

\begin{flalign}
\label{Qdef}
A_{{j}} &= n_x(j) \ \sin(\frac{\theta(j)}{2}), 
\\
B_{{j}} &= n_y(j) \ \sin(\frac{\theta(j)}{2}),
\\
C_{{j}} &= n_z(j) \ \sin(\frac{\theta(j)}{2}), 
\\
\label{Qdefend}
D_{{j}} &= \cos(\frac{\theta(j)}{2})
\end{flalign}

using the definitions of Eqs$.$ (\ref{eqtheta})-(\ref{eqthetaend}). Quaternions have the advantage that a series of rotations can be directly evaluated via the product

\begin{equation} {Q_2} \cdot {Q}_{1}=
\begin{pmatrix}
+{D_2} & -{C_2} & +{B_2} & +A_2 \\ 
+{C_2} & +{B_2} & -A_2 & +{B_2} \\
-{B_2} & +A_2 & +{D_2} & +{C_2} \\
-A_2 & -{B_2} & -{C_2} & +{D_2}
\end{pmatrix}\begin{pmatrix}
{A}_{1} \\ {B}_{1} \\ {C}_{1} \\ {D}_{1}
\end{pmatrix}\end{equation}

which, with its 16 simple multiplications and sums, is computationally more efficient than the construction of a rotation matrix involving at least one sine/cosine calculation. 

Using the formalism for universal rotation pulses, the cost function can directly be written as

\begin{equation}
\Phi_{\mathrm{UR}}=Re\braket{U_{\mathrm{F}}}{U_N}=Re\braket{Q_{\mathrm{F}}}{Q_N \cdots Q_j \cdots  Q_1} 
\end{equation}

and the gradients are derived as

\begin{equation}
\label{eq:gradient_urQ}
\Gamma_{\mathrm{UR}}(j)=\begin{pmatrix}\dfrac{\partial\Phi_{\mathrm{UR}}}{\partial\theta_x(j)} \\[4mm] \dfrac{\partial\Phi_{\mathrm{UR}}}{\partial\theta_y(j)} \\[4mm]
\dfrac{\partial\Phi_{\mathrm{UR}}}{\partial\theta_z(j)}
\end{pmatrix}
= \begin{pmatrix}Re \braket{Q_{\mathrm{F}}}{Q_N \cdots Q_{j+1} \cdot  
\dfrac{\partial Q_j}{\partial\theta_x(j)} \cdot Q_{j-1} \cdots Q_1 } \\[4mm]
Re \braket{Q_{\mathrm{F}}}{Q_N \cdots Q_{j+1} \cdot  
\dfrac{\partial Q_j}{\partial\theta_y(j)} \cdot Q_{j-1} \cdots Q_1 } \\[4mm]
Re \braket{Q_{\mathrm{F}}}{Q_N \cdots Q_{j+1} \cdot  
\dfrac{\partial Q_j}{\partial\theta_z(j)} \cdot Q_{j-1} \cdots Q_1 } 
\end{pmatrix},
\end{equation}

where again the overall gradient calculation is reduced to simple quaternion propagation plus the calculation of the rotation derivative - this time in quaternion notation. With Eqs$.$ (\ref{Qdef})-(\ref{Qdefend}) the involved partial derivatives of the individual quaternion components can be calculated and corresponding results are summarized in Table~\ref{tab:ur}. Obviously, the components of interest can also be derived with respect to polar coordinates $\alpha$ and $\theta_{xy}$ in the $xy$-plane, for which the terms are given in the right column of the same Table. It should be noted that the derivatives with respect to $\alpha$ are particularly simple, which guarantees fastest gradient calculation times for the case of constant amplitude pulses.

\subsection{Limited rf-amplitudes as holonomic constraints}

Optimizing efficient pulses in magnetic resonance usually implies optimizations with boundary conditions. One of the fundamental boundaries concern large rf-amplitudes that can cause experimental problems in the form of arching, violated duty cycles, and exceedance of allowed energy depositions. As such, they pose clear restrictions to applicable rf-amplitudes, rf-power, or rf-energy, respectively. Simple implementations with hard cut-off limits have been proposed early on in optimal control optimizations \cite{gershenzonLinearPhaseSlope2008, gershenzonOptimalControlDesign2007, kobzarExploringLimitsBroadband2004, kobzarExploringLimitsBroadband2008a, skinnerReducingDurationBroadband2004}. However, mathematically more sound are so-called holonomic constraints that can be included in optimizations as Lagrange multipliers. Holonomic constraints essentially require continuously differentiable functions for restrictions and are usually implemented by auxiliary variables. A multitude of periodic functions have been proposed as auxiliary functions \cite{goodwinAdvancedOptimalControl2017}, but we will concentrate here on a widely used trigonometric function, the hyperbolic tangent $\tanh$, which can be incorporated as a reduced rf-amplitude in polar coordinate systems according to

\begin{equation}
\theta_{xy}^{\rm{red}}(j) = \theta^{\rm{max}}(j) \ \cdot \ \tanh \left( \frac{\theta_{xy}(j)}{\theta^{\rm{max}}(j)} \right)
\end{equation} 

where the reduced amplitude $\theta_{xy}^{\rm{red}}(j)$ of the $j^{th}$ digit continuously converges only to the maximum value $\theta^{\rm{max}}(j)$. In an optimization, $\theta_{xy}(j)$ instead can now be changed as an auxiliary variable without any restriction, while $\theta_{xy}^{\rm{red}}(j)$ will be used for the actual values of rf-amplitudes in the pulse shape. 

For a simple amplitude restriction to a $j$-dependent maximum value, 

\begin{equation}
\theta^{\rm{max}}(j) = \mathrm{const} (j) = 2 \pi \ \Delta t (j) \ \nu_{\rm{rf}}^{\rm{max}}(j)
\end{equation}

with maximum rf-amplitude $\nu_{\rm{rf}}^{\rm{max}}(j)$ can be imposed. Usually a constant restriction for all time steps is chosen, but in special cases, physics may imply a $j$-dependent upper limit function \cite{skinnerNewStrategiesDesigning2012}. In an actual optimization, we restrict ourselves to the polar case where the rf-amplitude is inherently included in the variable $\theta_{xy}$. However, using the auxiliary approach, the variables $n_x$, $n_y$, $n_z$, $\theta_{xy}$, and $\theta$ from equations (\ref{eq16},\ref{eq17},\ref{eq24},\ref{eq25},\ref{eq27} - \ref{Qdefend}) used in the rotation matrices as well as in the quaternion formalism, have to be replaced by their reduced variants according to

\begin{eqnarray}
\theta_{xy}^{\rm{red}}(j) &=& \theta^{\rm{max}}(j) \ \cdot \ \tanh \left( \frac{\theta_{xy}(j)}{\theta^{\rm{max}}(j)} \right) \\
\theta^{\rm{red}}(j) &=& \sqrt{(\theta_{xy}^{\rm{red}}(j))^2 + \theta_z(j)^2} \\
n_x^{\rm{red}}(j) &=& \frac{\cos(\alpha) \cdot \theta_{xy}^{\rm{red}}(j)}{\theta^{\rm{red}}(j)} \\
n_y^{\rm{red}}(j) &=& \frac{\sin(\alpha) \cdot \theta_{xy}^{\rm{red}}(j)}{\theta^{\rm{red}}(j)} \\
n_x^{\rm{red}}(j) &=& \frac{\theta_{z}(j)}{\theta^{\rm{red}}(j)} 
\end{eqnarray}

The resulting gradients $\Gamma_{\rm{PP}}$ and $\Gamma_{\rm{UR}}$ are then essentially of the same form as derived above. Only for the derivation to the auxiliary variable itself the $\tanh$ is not a constant and care has to be taken. Taking the derivation of the $j^{th}$ element of an UR optimization as an example, we can derive

\begin{equation}
\frac{\partial \Phi_{\rm{UR}}}{\partial \theta_{xy} (j)} = \underbrace{
\frac{\partial \Phi_{\rm{UR}}}{\partial \theta_{xy}^{\rm{red}} (j)}}_{\textrm{control}} \, \cdot \, \underbrace{
\frac{d \theta_{xy}^{\rm{red}} (j)}{d \theta_{xy} (j)}}_{\textrm{auxiliary}}
\end{equation}

With the auxiliary derivative term the actual control $\theta_{xy}^{\rm{red}} (j)$ (i.e. the rf-amplitudes of the shaped pulse) can be used, for the optimization, instead, the variable $\theta_{xy}$ has to be used, making it necessary to include the terms

\begin{equation}
\frac{d \theta_{xy}^{\rm{red}} (j)}{d \theta_{xy} (j)} = \sech^2 \left( \frac{\theta_{xy}(j)}{\theta^{\rm{max}}(j)} \right) 
= 1 - \tanh^2 \left( \frac{\theta_{xy}(j)}{\theta^{\rm{max}}(j)} \right)
\end{equation}

into the equations of the gradient in question. Please note that the second solution is derived from the general relation $\sech^2(x) + \tanh^2(x) = 1$ and results in a computationally friendly term, as  $\tanh\left( \frac{\theta_{xy}(j)}{\theta^{\rm{max}}(j)} \right)$ has to be calculated anyway. 

If overall rf-power restrictions have to be imposed, the maximum rotation angle can also be  defined via the reduced angle at time point $j$ according to

\begin{equation}
\label{eqP1}
\theta_{xy}^{\rm{red}}(j) = \theta_{xy}(j) \ \cdot \ \sqrt{\frac{\overline{P^{\rm{max}}}}{\overline{P}}} \ \cdot \ \tanh \left( \sqrt{\frac{\overline{P}}{\overline{P^{\rm{max}}}}}  \right)
\end{equation} 

with the squareroots of the maximum allowed average rf-power $\overline{P^{\rm{max}}}$ and actual average rf-power expressed by

\begin{equation}
\overline{P} = \sum_{j = 1}^N \frac{(\theta_{xy}(j))^2}{N}
\end{equation}

where piecewise constant rf-amplitudes and uniform time steps $\Delta t$ have been assumed. Again, the derivative of the reduced rotation angle has to be calculated for the gradient, but this time the reduction involves rotation angles at all time steps and the derivative in e.g. an UR optimization with rf-power restriction will involve

\begin{equation}
\frac{\partial \Phi_{\rm{UR}}}{\partial \theta_{xy} (j)} = 
\sum_{k = 1}^N \ \frac{\partial \Phi_{\rm{UR}}}{\partial \theta_{xy}^{\rm{red}} (k)} \, \cdot \, 
\frac{\partial \theta_{xy}^{\rm{red}} (k)}{\partial \theta_{xy} (j)}
\end{equation}

with the reduced rotation angle derivatives defined by

\begin{eqnarray}
\frac{\partial \theta_{xy}^{\rm{red}}(j)}{\partial  \theta_{xy} (j)} & = &
\left( 1 - \frac{(\theta_{xy}(j))^2}{\sum_{i = 1}^N (\theta_{xy}(j))^2} \right)  \cdot  \sqrt{\frac{\overline{P^{\rm{max}}}}{\overline{P}}} \tanh \left( \sqrt{\frac{\overline{P}}{\overline{P^{\rm{max}}}}}  \right)  \nonumber \\
& &  + \frac{(\theta_{xy}(j))^2}{\sum_{i = 1}^N (\theta_{xy}(j))^2} \left( 1 - \tanh^2 \left( \sqrt{\frac{\overline{P}}{\overline{P^{\rm{max}}}}}  \right) \right) \ \ \  \ \forall \ \ \ k=j \ ;\\
\frac{\partial \theta_{xy}^{\rm{red}}(k)}{\partial  \theta_{xy} (j)} & = &
\frac{(\theta_{xy}(j))^2}{\sum_{i = 1}^N (\theta_{xy}(j))^2}   \cdot  \sqrt{\frac{\overline{P^{\rm{max}}}}{\overline{P}}} \tanh \left( \sqrt{\frac{\overline{P}}{\overline{P^{\rm{max}}}}}  \right)  \nonumber \\
& &  + \frac{(\theta_{xy}(j))^2}{\sum_{i = 1}^N (\theta_{xy}(j))^2} \left( 1 - \tanh^2 \left( \sqrt{\frac{\overline{P}}{\overline{P^{\rm{max}}}}}  \right) \right) \ \ \  \ \forall  \ \ \ k \ne j \ . 
\label{eqP2}
\end{eqnarray} 

Finally, also the overall rf-energy of a pulse can be restricted, using, for example, the energy defined in Hz

\begin{equation}
\frac{E}{h} = \overline{P} \, t_p,
\end{equation}

or the more convenient expression in terms of rotation angles

\begin{equation}
E_{\theta} = \sum_{i = 1}^N (\theta_{xy}(i))^2 = (2 \pi)^2 \ \Delta t \ \frac{E}{h}
\end{equation}

when equations (\ref{eqP1} - \ref{eqP2}) are used with substituting $\overline{P^{\rm{max}}}$ with the maximum allowed energy $E_{\theta}^{\rm{max}} =  (2 \pi)^2 \Delta t E^{\rm{max}} / h$ and  $\overline{P}$ by $E_{\theta}$.

\section{Results and Discussion}

\subsection{Runtime comparisons}

All analytical solutions for the different gradients have been derived to significantly enhance computational performance. We therefore implemented the analytical solutions in Julia, a modern compiling programming language with a comfortable interface similar to Python or Matlab. The implementation of linear algebra routines in Julia is furthermore based on efficient BLAS routines, that still represent today's standard. We compared the analytical solutions derived in the various Tables with the exact solution provided by the augmented matrix exponential \cite{floetherRobustQuantumGates2012, goodwinAuxiliaryMatrixFormalism2015, vanloanComputingIntegralsInvolving1978} approach and, in addition, to a finite difference approximation. All calculations were performed with double precision (Float64) accuracy. 

For the Cartesian point-to-point case, the augmented matrix exponential approach for a specific Cartesian control can be implemented in a 6-dimensional matrix that consists of two identical rotation matrices along the diagonal, and the generator of corresponding the $x$, $y$, or $z$-control in the right upper corner. Solving the equation

\begin{equation}
\begin{pmatrix} 
R_j & \frac{\partial R_j}{\partial \theta_{\alpha}(j)} \\
0 & R_j
\end{pmatrix}
= \exp 
\begin{pmatrix} 
-i L \theta(j)  & -i G_{\alpha} \Delta t \\
0 &  -i L \theta(j)
\end{pmatrix}
\ \ \ \rm{with} \ \ \   L \theta(j) = (G_x \theta_x(j) + G_y \theta_y(j) +G_z \theta_z(j)) 
\ \ \ , \ \ \ \alpha \in \{x, y, z\}
\end{equation}

and 
\begin{equation}
G_x = \begin{pmatrix} 
0 & 0 & 0 \\
0 & 0 & 1 \\
0 & 1 & 0
\end{pmatrix}
\ \ \ , \ \ \
G_y = \begin{pmatrix} 
0 & 0 & 1 \\
0 & 0 & 0 \\
1 & 0 & 0
\end{pmatrix}
\ \ \ , \ \ \
G_z = \begin{pmatrix} 
0 & 1 & 0 \\
1 & 0 & 0 \\
0 & 0 & 0
\end{pmatrix}
\end{equation}

 gives the desired result for the first derivative $\partial R_j / \partial \theta_{\alpha}$ in place of the control. 
It should be noted, that the augmented matrix exponential approach cannot directly be applied to polar coordinates, as no generator is available for a change of the phase-angle $\alpha$ or the overall rotation angle $\theta_{{xy}}$. This would imply the calculation of an x- and y-control and then combine the two results for obtaining phase and amplitude. We did not attempt this, as this clearly would result in doubling of calculation time and an even larger gap in performance of this approach. Correspondingly, the auxiliary $\tanh$ approach was not attempted, as it only makes sense in polar coordinates. 

Equivalent to rotational matrices, the augmented matrix exponential approach can be applied using a single spin propagator defined via the Pauli matrices. The corresponding equation involves $4 \times 4$ matrices according to 

\begin{equation}
\begin{pmatrix} 
U_j & \frac{\partial U_j}{\partial \theta_{\alpha}(j)} \\
0 & U_j
\end{pmatrix}
= \exp 
\begin{pmatrix} 
-i H \theta(j)  & -i \sigma_{\alpha} \Delta t \\
0 &  -i H \theta(j)
\end{pmatrix}
\ \ \ \rm{with} \ \ \ \sigma \theta(j) = (\sigma_x \theta_x(j) + \sigma_y \theta_y(j) +\sigma_z \theta_z(j)) 
\ \ \ , \ \ \ \alpha \in \{x, y, z\}
\end{equation}

and 
\begin{equation}
\sigma_x = \begin{pmatrix} 
 0 & \frac{1}{2} \\
 \frac{1}{2} & 0
\end{pmatrix}
\ \ \ , \ \ \
\sigma_y = \begin{pmatrix} 
 0 & -\frac{i}{2} \\
 \frac{i}{2} & 0
\end{pmatrix}
\ \ \ , \ \ \
\sigma_z = \begin{pmatrix} 
 \frac{1}{2} & 0 \\
 0 & -\frac{1}{2}
\end{pmatrix}.
\end{equation}

Finite difference approximations were performed using 3D rotation matrices and quaternions in the PP and UR cases, respectively. They essentially consist of two rotations and the corresponding calculations are mainly determined by the computation of trigonometric terms. Averaged computation times were calculated using the Julia package BenchmarkTools. The resulting times averaged over 1000 derivative calculations are summarized in Table \ref{tab:runt}. In all cases the analytical derivative solutions introduced here  outperforms the other approaches. Compared to the other exact calculations based on the augmented matrix exponential approach, improvements of factors 100-150 are observed. The finite differences approach, on the other hand, is quite fast with factors of approximately 1-3 for the different types of optimization. The gain of the analytical solution compared in this case is not so much in the speed, but, of course, in the higher accuracy of the gradient and therefore a faster overall convergence of the gradient-based optimizations.  

\subsection{Example optimizations}

To test the performance of the analytical exact gradients, we also looked into three different optimization scenarios. The first one concerns with a relatively simple, but nevertheless important application: the development of  $^{15}$N pulses that cover the full 50~ppm amide bandwidth on a 1.2~GHz spectrometer with a maximum rf-amplitude that corresponds to a rectangular 90$^{\circ}$ pulse of 50~$\mu$s. The aim was not to produce the best pulse possible, but to monitor and compare how well the analytical solutions were able to produce a good pulse and what the average duration per iteration as a measure of computation speed was. The results are summarized in Table \ref{tab:15N} for excitation (so-called BEBOP pulses \cite{skinnerApplicationOptimalControl2003b, skinnerTailoringOptimalControl2005a,kobzarExploringLimitsBroadband2012a}) and inversion (so-called BIBOP pulses \cite{ehniBEBEtrBUBIJcompensated2013b, kobzarExploringLimitsBroadband2004, kobzarExploringLimitsBroadband2012a, nimbalkarFantasticFourPlug2013}) pulse optimizations. 

\begin{figure}[t]
\includegraphics[width=8.5cm]{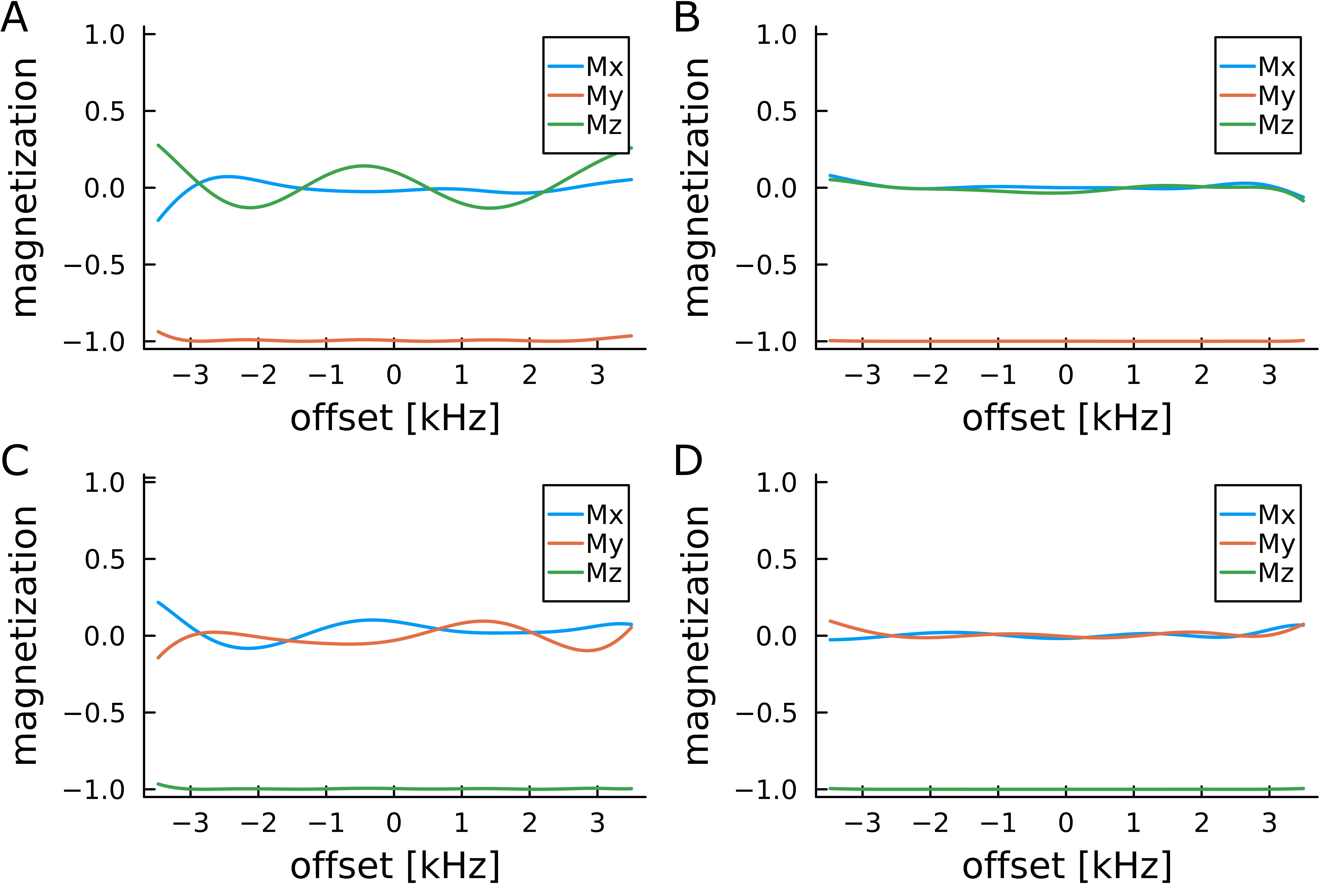}
\caption{Offset profiles for excitation (A,B) and inversion (C,D) pulses optimized for amide nitrogen excitation of proteins on a 1.2~GHz spectrometer. The quality factors of the corresponding pulse shapes are written in italics in Table \ref{tab:15N}. For both types, worst (A,C) and best (B,D) pulses were simulated and differences are clearly visible.
\label{fig1}
}
\end{figure}

For parameters we decided to focus on amplitude-restricted pulses with a maximum rf-amplitude of 5~kHz and $\pm$10\% B$_1$-compensation. The bandwidth is 6000~kHz, which is only slightly larger than the rf-amplitude of the pulse. For the optimization we used the Optim.jl Julia package (https://julianlsolvers.github.io/Optim.jl/stable/), with its BFGS implementation. We ran all optimizations on the very same laptop with which we compared the runtime performance in Table \ref{tab:runt}. Independent of the parameters used for optimization, the resulting pulses generally gave high quality factors larger than 0.99. Three different digitizations were used in the optimizations, where durations of the digits of 1, 10, and 50~$\mu$s led to 500, 50, and 10 piecewise constant pulse elements, respectively. It is no surprise that the latter led to fastest overall optimization times in the sub-second range, as well as the shortest computation times per iteration and also on average less iterations for convergence. Interestingly, the resulting quality factors for the short optimizations are overall in the same range as for the other optimizations. Particularly for initial screens it might therefore be of use to start a large number of fast optimizations with long $\Delta t$ elements, which then form a good basis to choose optimization parameters for a more detailed search. However, it should be noted, that there are profound differences between pulses with quality factors on the order of 0.992 and ones with quality factors of 0.999. This can nicely be seen in Fig. \ref{fig1}, where corresponding excitation and inversion pulses are compared in their performance. Clearly the higher quality factors correspond to significantly better performing pulse shapes. 

Comparing the different controls, it is no surprise that constant amplitude optimizations that only control the pulse phase $\alpha$, resulted in shortest calculation times. This approach simply has the least amount of controls and also the smallest parameter space in the optimizations. More unexpectedly, the second fastest overall optimization times were reached by the reduced amplitude parameters $\theta_{xy}^{\rm{red}}$ using the tanh-based restriction of equation (\ref{eqP1}). Although the additional computation of the tanh terms is computationally costly, it seems to clearly outperform the Lagrange multiplier type amplitude restriction used in the other cases. It seems therefore worth following the pseudo parameter approach when dealing with such type of restrictions. 

For the more conventional approaches using either $x$, $y$ or $\theta_{xy}$, $\alpha$ controls, no real difference in performance can be seen. Also the addition of $z$-controls did not really change the convergence characteristic. However, we experienced that $z$-controls lead to many optimizations that stop at a very early stage of an optimization. Apparently, the unrestricted $z$-controls lead sometimes to erratic behaviour with large frequency jumps between neighbouring digits. In principle, $z$-control should allow better convergence and a larger parameter space when rough digitizations are used in pulse optimization. The situation might change, when additional constraints like smoothness restrictions are applied to $z$-controls.   

\begin{figure}[t]
\includegraphics[width=8.5cm]{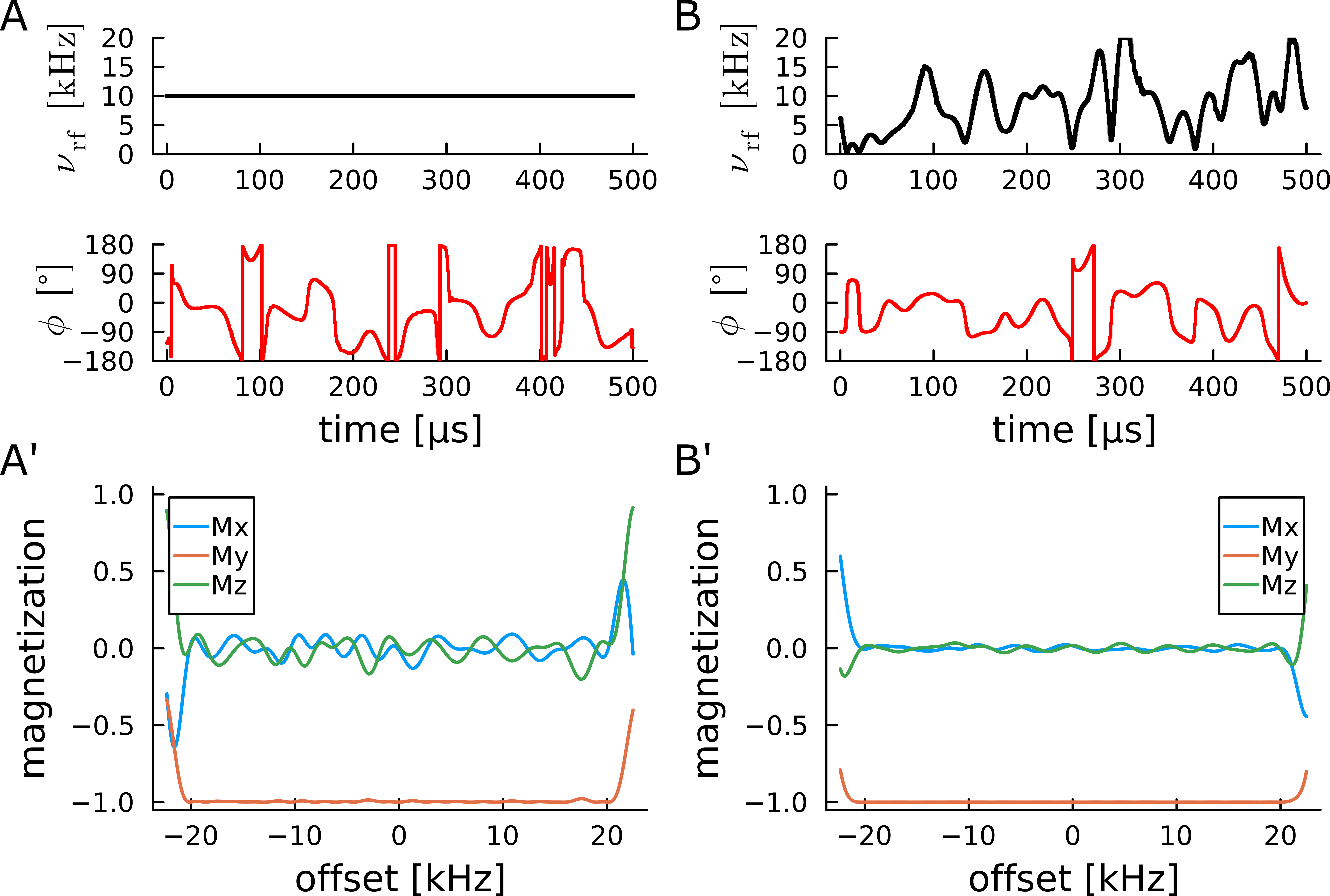}
\caption{Offset profiles for two excitation pulses with identical overall rf-power (A,B) and optimization parameters given in Table \ref{tab:13C}. The best obtained BEBOP pulse shape with constant amplitude and overall quality factor 0.99715 has significantly lower performance (A') than the power-BEBOP with quality factor 0.99982 (B'), corroborating well-known results regarding physical limits studies \cite{kobzarExploringLimitsBroadband2008a}. 
\label{fig2}
}
\end{figure}

\begin{figure}[t]
\includegraphics[width=8.5cm]{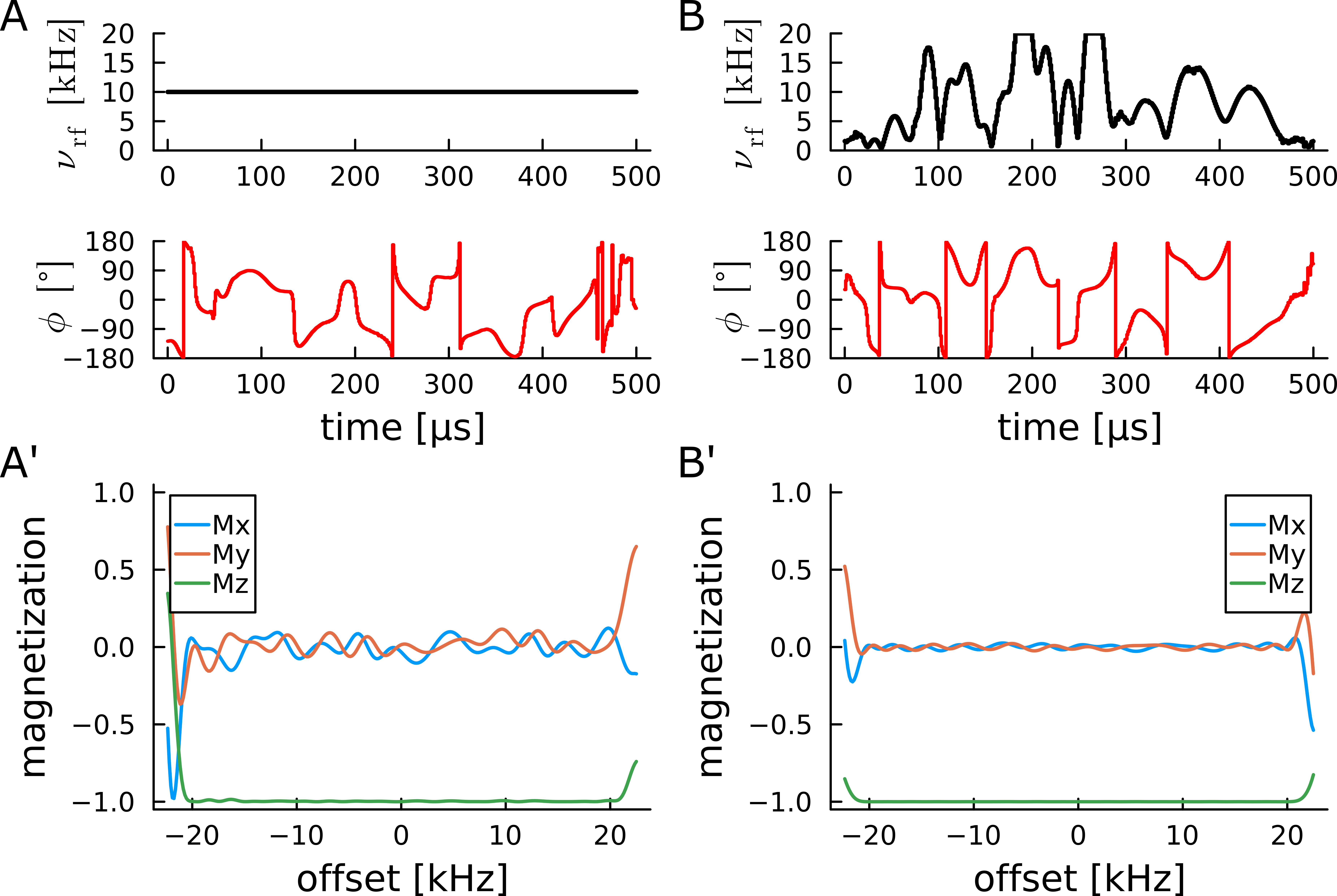}
\caption{Offset profiles for two inversion pulses with identical overall rf-power (A,B) and optimization parameters given in Table \ref{tab:13C}. The best obtained BIBOP pulse shape with constant amplitude and overall quality factor 0.99744 has significantly lower performance (A') than the power-BIBOP with quality factor 0.99983 (B'), corroborating well-known results regarding physical limits studies \cite{kobzarExploringLimitsBroadband2008a}.
\label{fig3}
}
\end{figure}

\begin{figure}[t]
\includegraphics[width=8.5cm]{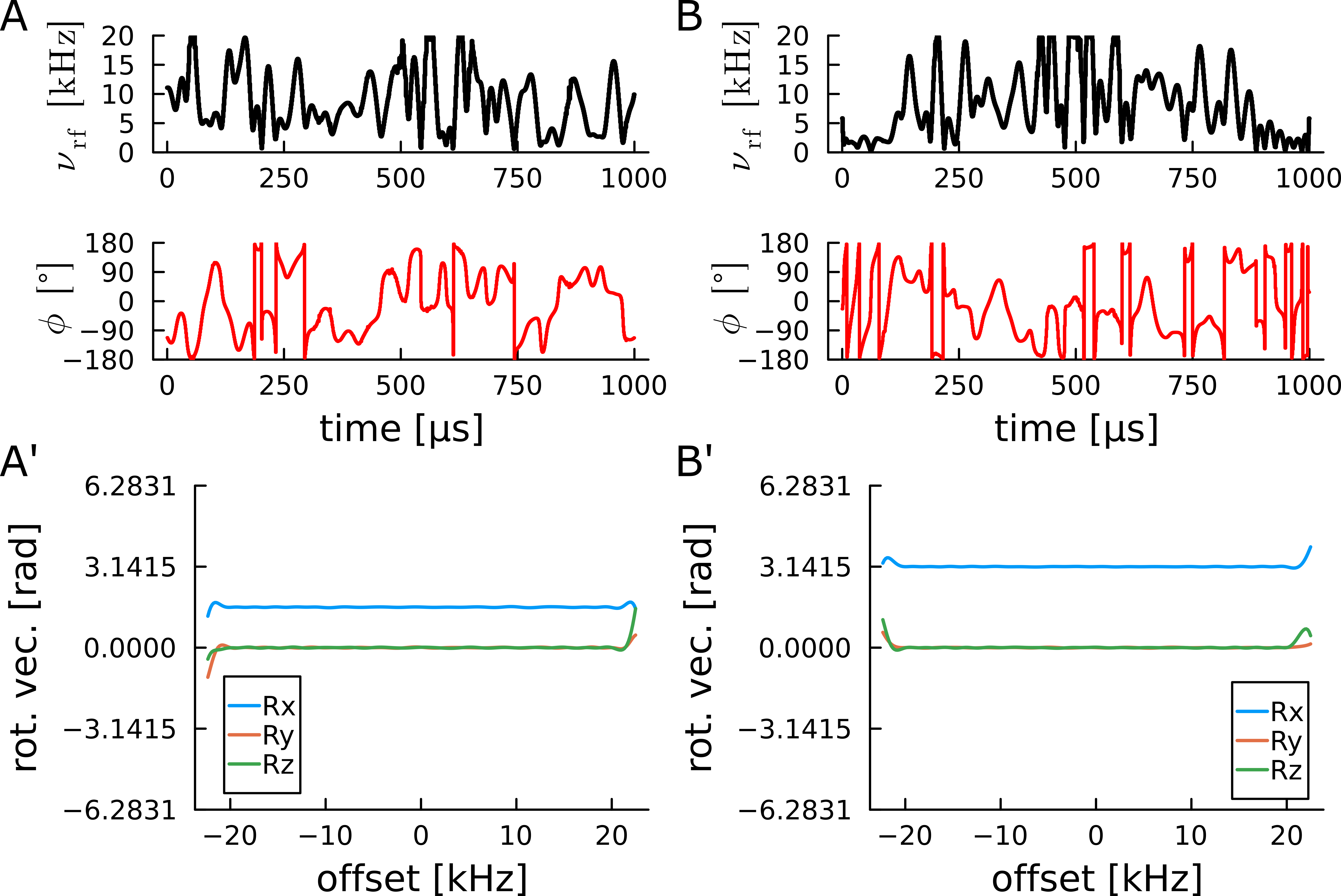}
\caption{Offset profiles for a universal rotation power-BURBOP-90$^{\circ}$ (A) and a universal rotation power-BURBOP-180$^{\circ}$ (B) pulse with identical overall rf-power and optimization parameters given in Table \ref{tab:13C}. Both pulses show exceptional performance.
\label{fig4}
}
\end{figure}

The second studied scenario concerns pulses with a bandwidth of 40~kHz, which corresponds to a 266~ppm carbon chemical shift range on a 600~MHz spectrometer, and to 160~ppm on a 1.0~GHz spectrometer. We optimized a basic set of excitation, inversion and different universal rotation pulses with either constant amplitude of 10~kHz, or with a combined Lagrange multiplier type rf-amplitude (20~kHz) and rf-power ($sqrt{\overline{P^{\rm{max}}}}$) restrictions. Corresponding pulse shapes generally resulted in very good performance pulses, as can also be seen in Figs. \ref{fig2}-\ref{fig4} for the different pulse types. In all cases the amplitude and power-restricted pulse shapes perform better than the constant amplitude pulses with identical overall rf-power, which could be expected from the results of previous systematic studies of BEBOP/BIBOP \cite{kobzarExploringLimitsBroadband2004} and power-BEBOP/power-BIBOP \cite{kobzarExploringLimitsBroadband2008a} pulses. 

\begin{figure*}[t]
\includegraphics[width=19cm]{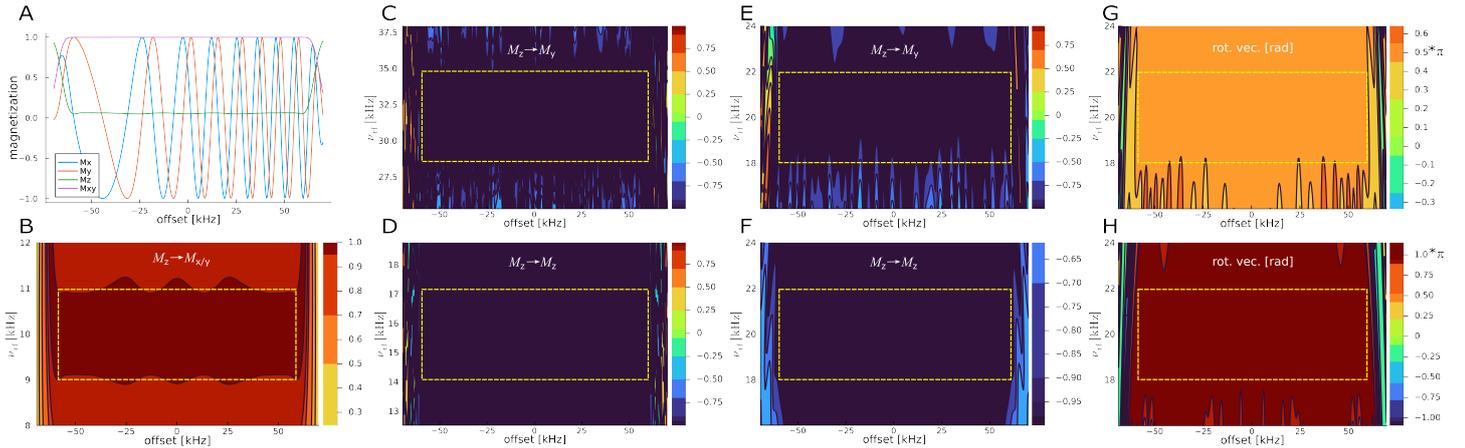}
\caption{Summary of pulse shapes optimized for fluorine excitation. (A,B) extremely short saturation pulse of 120~$\mu$s covering the full bandwidth of 120 kHz. Although the pulse was optimized without B$_1$-compensation, considerable compensation is achieved in the $\pm$10\% range. (C,D) power-BEBOP excitation and power-BIBOP inversion pulses with the squareroot of the average rf-power of only 10~kHz and 1.4~ms duration. (E,F) corresponding BEBOP excitation and BIBOP inversion pulses with constant rf-amplitude of 20~kHz and durations of 350~$\mu$s and below. (G,H) Universal rotation BURBOP-90 and BURBOP-180 pulses of 600~$\mu$s and 700~$\mu$s duration equivalent to previously published pulse shapes \cite{lingelComprehensiveHighThroughputExploration2020}. The Offset/B$_1$ region used for most pulse optimizations is highlighted by the yellow dashed boxes. 
\label{fig5}
}
\end{figure*}

Finally, we optimized a set of pulses for $^{19}$F spectroscopy on a 600~MHz spectrometer, where previously two universal rotation pulses were optimized for screening experiments with essentially identical parameters \cite{lingelComprehensiveHighThroughputExploration2020}. The large bandwidth and considerable B$_1$-compensation for most pulses is remarkable and nicely documented in Fig.~\ref{fig5}. 

A first pulse is a saturation pulse bringing $z$ magnetization efficiently into the $x$,$y$ plane. As has been previously experienced \cite{endersCytotoxicityNMRStudies2014, luyUltrabroadbandNMRSpectroscopy2009}, the wide offset range can be covered with a very short pulse shape. It should also be noted that the overall optimization time of this difficult pulse on a single processor core required only 1.6~seconds with a close to perfect performance. Although no B$_1$-compensation was applied, the pulse will be very well applicable, as can be seen in Fig.~\ref{fig5}B. Surprisingly good are also the excitation and inversion pulses with their overall rf-power restricted to the equivalent to a 10~kHz constant amplitude pulse. The restriction in this case has been applied via the tanh-type pseudo controls with inherent power-based scaling. Both pulse shapes with a high number of offset checks, 700 piecewise constant pulse digits, and detailed B$_1$-compensation were the most complex optimizations attempted here. Correspondingly, optimization times comprise half a day to a day on a single core. The bandwidth-over-B$_1$-ratio in this case is 12, which is much larger than for most previously optimized excitation and inversion pulses, and demonstrates to some extent the increase in complexity and optimization duration with more demanding pulse shapes. Both optimizations did not fully converge, but stopped after the maximum number of iterations specified (5000 iterations). It can be expected that optimizations that run until convergence is reached will lead to better pulse shapes, but most likely with only slight improvements. Finally, constant amplitude pulse shapes with a bandwidth-over-B$_1$-ratio of 6 were optimized for excitation, inversion, and universal rotation 90$^{\circ}$ and 180$^{\circ}$ pulses. All of the pulses are rather short for their performance and should be readily suited for corresponding $^{19}$F-based applications. 

\section{Conclusion}

Cost functions and corresponding exact gradients for point-to-point and universal rotation type single spin optimizations have been derived analytically. Their implementation into the self-written program PulseOptimizer (will be published somewhere else) resulted into improvements in computation times compared to other analytical gradient calculations like the augmented matrix exponentiation approach by roughly two orders of magnitude. Correspondingly, very different types of conventional pulse shapes could be optimized in a matter of seconds to few hours using a single processor core on a laptop. With this increase in optimization speed an important step towards the optimization of pulse shapes with very large bandwidth-over-B$_1$-ratios is done. However, already some of the example pulses have shown that further technical improvements in optimization software will be necessary to tackle pulses with bandwidths that are 100 times larger than the corresponding B$_1$ field strenghts. As such, massive parallelization and eventually efficient implementation in GPU architectures either on a quality factor level and/or on an optimizer level seem to be appropriate additional future steps.    

\section{Acknowledgement}

B.L. thanks the Deutsche Forschungsgemeinschaft (CRC 1527 - HyPERION, project C01) and the HGF programme MSE Information (43.35.02) for financial support. We also thank very much David L. Goodwin (Aarhus and Oxford) for many detailed and decisive discussions. 

\bibliographystyle{unsrt}      
\bibliography{test2}   

\begin{table}[hbt!]
\caption {Derivatives of rotation matrix components with respect to $\theta_{x}$ and $\theta_{y}$ for point-to-point optimizations.}
\label{tab:pp_dx} 
\begin{tabular}{| rl |}
  \hline 
& \\
 \ \ $\dfrac{\partial R_{xx}}{\partial \theta_x}$ & $= (n_x^3 - n_x)\Bigg( \sin(\theta) + \dfrac{2 \ (\cos(\theta)-1)}{\theta} \Bigg)$ \\[5 mm]
$\dfrac{\partial R_{xy}}{\partial \theta_x}$ & $=\dfrac{n_y-2 n_x^2n_y}{\theta} + \Bigg( \dfrac{2 n_x^2n_y-n_y}{\theta} - n_xn_z \Bigg) \cos(\theta) + \Bigg( n_x^2n_y +\dfrac{n_xn_z}{\theta} \Bigg) \sin(\theta) $ \\[5mm]
$\dfrac{\partial R_{xz}}{\partial \theta_x}$ & $=\dfrac{n_z-2 n_x^2n_z}{\theta} + \Bigg( \dfrac{2 n_x^2n_z-n_z}{\theta} + n_xn_y \Bigg) \cos(\theta) + \Bigg( n_x^2n_z -\dfrac{n_xn_y}{\theta} \Bigg) \sin(\theta) $ \\[5mm]
$\dfrac{\partial R_{yx}}{\partial \theta_x}$ & $=\dfrac{n_y-2 n_x^2n_y}{\theta} + \Bigg( \dfrac{2 n_x^2n_y - n_y}{\theta} + n_xn_z \Bigg) \cos(\theta) + \Bigg( n_x^2n_y -\dfrac{n_xn_z}{\theta} \Bigg) \sin(\theta) $ \\[5mm]
$\dfrac{\partial R_{yy}}{\partial \theta_x}$ & $=(n_xn_y^2-n_x)\sin(\theta)+n_xn_y^2 \ \Bigg(\dfrac{2 \ \left(\cos(\theta)-1\right)}{\theta}\Bigg)$
\\[5 mm]
$\dfrac{\partial R_{yz}}{\partial \theta_x}$ & $=- \dfrac{2 n_x n_y n_z}{\theta} + \Bigg( \dfrac{2 n_xn_yn_z}{\theta} - n_x^2 \Bigg) \cos(\theta) + \Bigg( n_xn_yn_z -\dfrac{(n_y^2+n_z^2)}{\theta} \Bigg) \sin(\theta) $ \\[5mm]
$\dfrac{\partial R_{zx}}{\partial \theta_x}$ & $=\dfrac{n_z-2 n_x^2n_z}{\theta} + \Bigg( \dfrac{2 n_x^2n_z-n_z}{\theta} - n_xn_y \Bigg) \cos(\theta) + \Bigg( n_x^2n_z +\dfrac{n_xn_y}{\theta} \Bigg) \sin(\theta) $ \\[5mm]
$\dfrac{\partial R_{zy}}{\partial \theta_x}$ & $=- \dfrac{2 n_x n_y n_z}{\theta} + \Bigg( \dfrac{2 n_xn_yn_z}{\theta} + n_x^2 \Bigg) \cos(\theta) + \Bigg( n_xn_yn_z +\dfrac{(n_y^2+n_z^2)}{\theta} \Bigg) \sin(\theta) $ \\[5mm]
$\dfrac{\partial R_{zz}}{\partial \theta_x}$ & $=(n_xn_z^2-n_x)\sin(\theta)+n_xn_z^2 \ \Bigg(\dfrac{2 \ \left(\cos(\theta)-1\right)}{\theta}\Bigg)$ \\[5mm]
\hline 
& \\
$\dfrac{\partial R_{xx}}{\partial \theta_y}$ & $=(n_yn_x^2-n_y)\sin(\theta)+n_yn_x^2 \ \Bigg(\dfrac{2 \ \left(\cos(\theta)-1\right)}{\theta}\Bigg)$ \\[5mm]
$\dfrac{\partial R_{xy}}{\partial \theta_y}$ & $=\dfrac{n_x-2 n_y^2n_x}{\theta} + \Bigg( \dfrac{2 n_y^2n_x-n_x}{\theta} - n_yn_z \Bigg) \cos(\theta) + \Bigg( n_y^2n_x +\dfrac{n_yn_z}{\theta} \Bigg) \sin(\theta) $ \\[5mm]
$\dfrac{\partial R_{xz}}{\partial \theta_y}$ & $=- \dfrac{2 n_x n_y n_z}{\theta} + \Bigg( \dfrac{2 n_xn_yn_z}{\theta} + n_y^2 \Bigg) \cos(\theta) + \Bigg( n_xn_yn_z +\dfrac{(n_x^2+n_z^2)}{\theta} \Bigg) \sin(\theta) $ \\[5mm]
$\dfrac{\partial R_{yx}}{\partial \theta_y}$ & $=\dfrac{n_x-2 n_y^2n_x}{\theta} + \Bigg( \dfrac{2 n_y^2n_x-n_x}{\theta} + n_yn_z \Bigg) \cos(\theta) + \Bigg( n_y^2n_x -\dfrac{n_yn_z}{\theta} \Bigg) \sin(\theta) $ \\[5mm]
 \ \ $\dfrac{\partial R_{yy}}{\partial \theta_y}$ & $= (n_y^3 - n_y)\Bigg( \sin(\theta) + \dfrac{2 \ (\cos(\theta)-1)}{\theta} \Bigg)$ \\[5 mm]
$\dfrac{\partial R_{yz}}{\partial \theta_y}$ & $=\dfrac{n_z-2 n_y^2n_z}{\theta} + \Bigg( \dfrac{2 n_y^2n_z-n_z}{\theta} - n_xn_y \Bigg) \cos(\theta) + \Bigg( n_y^2n_z +\dfrac{n_xn_y}{\theta} \Bigg) \sin(\theta) $ \\[5mm]
$\dfrac{\partial R_{zx}}{\partial \theta_y}$ & $=- \dfrac{2 n_x n_y n_z}{\theta} + \Bigg( \dfrac{2 n_xn_yn_z}{\theta} - n_y^2 \Bigg) \cos(\theta) + \Bigg( n_xn_yn_z -\dfrac{(n_x^2+n_z^2)}{\theta} \Bigg) \sin(\theta) $ \\[5mm]
$\dfrac{\partial R_{zy}}{\partial \theta_y}$ & $=\dfrac{n_z-2 n_y^2n_z}{\theta} + \Bigg( \dfrac{2 n_y^2n_z-n_z}{\theta} + n_xn_y \Bigg) \cos(\theta) + \Bigg( n_y^2n_z -\dfrac{n_xn_y}{\theta} \Bigg) \sin(\theta) $ \\[5mm]
$\dfrac{\partial R_{zz}}{\partial \theta_y}$ & $=(n_yn_z^2-n_y)\sin(\theta)+n_yn_z^2 \ \Bigg(\dfrac{2 \ \left(\cos(\theta)-1\right)}{\theta}\Bigg)$ \\[5mm]
\hline
\end{tabular}
\end{table}


\begin{table}[hbt!]
\caption {Derivatives of rotation matrix components with respect to $\theta_{z}$ for point-to-point optimizations.}
 \label{tab:pp_dz} 
\begin{tabular}{| rl |}
  \hline 
& \\
 \ \ $\dfrac{\partial R_{xx}}{\partial \theta_z}$ & $=(n_zn_x^2-n_z)\sin(\theta)+n_zn_x^2 \ \Bigg(\dfrac{2 \ \left(\cos(\theta)-1\right)}{\theta}\Bigg)$ \\[5mm]
$\dfrac{\partial R_{xy}}{\partial \theta_z}$ & $=- \dfrac{2 n_x n_y n_z}{\theta} + \Bigg( \dfrac{2 n_xn_yn_z}{\theta} - n_z^2 \Bigg) \cos(\theta) + \Bigg( n_xn_yn_z -\dfrac{(n_y^2+n_x^2)}{\theta} \Bigg) \sin(\theta) $ \\[5mm]
$\dfrac{\partial R_{xz}}{\partial \theta_z}$ & $=\dfrac{n_x-2 n_z^2n_x}{\theta} + \Bigg( \dfrac{2 n_z^2n_x-n_x}{\theta} + n_zn_y \Bigg) \cos(\theta) + \Bigg( n_z^2n_x -\dfrac{n_zn_y}{\theta} \Bigg) \sin(\theta) $ \\[5mm]
$\dfrac{\partial R_{yx}}{\partial \theta_z}$ & $=- \dfrac{2 n_x n_y n_z}{\theta} + \Bigg( \dfrac{2 n_xn_yn_z}{\theta} + n_z^2 \Bigg) \cos(\theta) + \Bigg( n_xn_yn_z +\dfrac{(n_y^2+n_x^2)}{\theta} \Bigg) \sin(\theta) $ \\[5mm]
$\dfrac{\partial R_{yy}}{\partial \theta_z}$ & $=(n_zn_y^2-n_z)\sin(\theta)+n_zn_y^2 \ \Bigg(\dfrac{2 \ \left(\cos(\theta)-1\right)}{\theta}\Bigg)$
\\[5 mm]
$\dfrac{\partial R_{yz}}{\partial \theta_z}$ & $=\dfrac{n_y-2 n_z^2n_y}{\theta} + \Bigg( \dfrac{2 n_z^2n_y - n_y}{\theta} - n_xn_z \Bigg) \cos(\theta) + \Bigg( n_z^2n_y +\dfrac{n_xn_z}{\theta} \Bigg) \sin(\theta) $ \\[5mm]
$\dfrac{\partial R_{zx}}{\partial \theta_z}$ & $=\dfrac{n_x-2 n_z^2n_x}{\theta} + \Bigg( \dfrac{2 n_z^2n_x-n_x}{\theta} - n_zn_y \Bigg) \cos(\theta) + \Bigg( n_z^2n_x +\dfrac{n_zn_y}{\theta} \Bigg) \sin(\theta) $ \\[5mm]
$\dfrac{\partial R_{zy}}{\partial \theta_z}$ & $=\dfrac{n_y-2 n_z^2n_y}{\theta} + \Bigg( \dfrac{2 n_z^2n_y - n_y}{\theta} + n_xn_z \Bigg) \cos(\theta) + \Bigg( n_z^2n_y -\dfrac{n_xn_z}{\theta} \Bigg) \sin(\theta) $ \\[5mm]
$\dfrac{\partial R_{zz}}{\partial \theta_z}$ & $= (n_z^3 - n_z)\Bigg( \sin(\theta) + \dfrac{2 \ (\cos(\theta)-1)}{\theta} \Bigg)$ \\[5 mm]
\hline
\end{tabular}
\end{table}

\begin{table}[hbt!]
\caption {Derivatives of rotation matrix components with respect to phase $\alpha$ and rf-amplitude $\theta_{xy}$  for point-to-point optimizations}
\label{tab:pp_awz1} 
\begin{tabular}{| rl |}
  \hline 
&  \\
 \ \ $\dfrac{\partial R_{xx}}{\partial \alpha}$&$=-2n_xn_y(1 - \cos(\theta))$ \\[5mm]
$\dfrac{\partial R_{xy}}{\partial \alpha}$&$=(n_x^2 - n_y^2)(1 - \cos(\theta))$ \\[5mm]
$\dfrac{\partial R_{xz}}{\partial \alpha}$&$=-n_yn_z(1 - \cos(\theta)) + n_x\sin(\theta)$ \\[5mm]
$\dfrac{\partial R_{yx}}{\partial \alpha}$&$=(n_x^2 - n_y^2)(1 - \cos(\theta))$ \\[5mm]
$\dfrac{\partial R_{yy}}{\partial \alpha}$&$=2n_xn_y(1 - \cos(\theta))$ \\[5mm]
$\dfrac{\partial R_{yz}}{\partial \alpha}$&$=n_xn_z(1 - \cos(\theta)) + n_y\sin(\theta)$ \\[5mm]
$\dfrac{\partial R_{zx}}{\partial \alpha}$&$=-n_yn_z(1 - \cos(\theta)) - n_x\sin(\theta)$ \\[5mm]
$\dfrac{\partial R_{zy}}{\partial \alpha}$&$=n_xn_z(1 - \cos(\theta)) - n_y\sin(\theta)$ \\[5mm]
$\dfrac{\partial R_{zz}}{\partial \alpha}$&$=0$ \\[3mm]
\hline
&  \\
$\dfrac{\partial R_{xx}}{\partial \theta_{xy}}$&$=-\dfrac{\theta_{xy}(n_y^2+n_z^2)}{\theta}\sin(\theta)+\dfrac{2n_x^2n_z^2}{\theta_{xy}}(1-\cos(\theta))$ \\[5mm]
$\dfrac{\partial R_{xy}}{\partial \theta_{xy}}$&$=-n_{xy} n_{z}\left(\cos(\theta) - \dfrac{\sin(\theta)}{\theta}\right) + n_{x} n_{y} n_{xy} \sin(\theta) + 2 \left(\dfrac{n_{x} n_{y}}{\theta_{xy}} - \dfrac{n_{x} n_{y} n_{xy}}{\theta}\right) (1 - \cos(\theta))$ \\[5mm]
$\dfrac{\partial R_{xz}}{\partial \theta_{xy}}$&$=\left(\dfrac{n_{y}}{\theta_{xy}} - \dfrac{n_{y} n_{xy}}{\theta} + n_{x} n_{z} n_{xy}\right) \sin(\theta) + n_{y} n_{xy} \cos(\theta) + \left(\dfrac{n_{x} n_{z}}{\theta_{xy}} - \dfrac{2 n_{x} n_{z} n_{xy}}{\theta}\right) (1 - \cos(\theta))$ \\[5mm]
$\dfrac{\partial R_{yx}}{\partial \theta_{xy}}$&$=n_{xy} n_{z}\left(\cos(\theta) - \dfrac{\sin(\theta)}{\theta}\right) + n_{x} n_{y} n_{xy} \sin(\theta) + 2 \left(\dfrac{n_{x} n_{y}}{t_{xy}} - \dfrac{n_{x} n_{y} n_{xy}}{\theta}\right) (1 - \cos(\theta))$ \\[5mm]
$\dfrac{\partial R_{yy}}{\partial \theta_{xy}}$&$=-\dfrac{\theta_{xy}(n_x^2+n_z^2)}{\theta}\sin(\theta)+\dfrac{2n_y^2n_z^2}{\theta_{xy}}(1-\cos(\theta))$ \\[5mm]
$\dfrac{\partial R_{yz}}{\partial \theta_{xy}}$&$=\left(-\dfrac{n_{x}}{\theta_{xy}} + \dfrac{n_{x} n_{xy}}{\theta} + n_{y} n_{z} n_{xy}\right) \sin(\theta) - n_{x} n_{xy} \cos(\theta) + \left(\dfrac{n_{y} n_{z}}{\theta_{xy}} - \dfrac{2 n_{y} n_{z} n_{xy}}{\theta}\right) (1 - \cos(\theta))$ \\[5mm]
$\dfrac{\partial R_{zx}}{\partial \theta_{xy}}$&$=\left(-\dfrac{n_{y}}{\theta_{xy}} + \dfrac{n_{y} n_{xy}}{\theta} + n_{x} n_{z} n_{xy}\right) \sin(\theta) - n_{y} n_{xy} \cos(\theta) + \left(\dfrac{n_{x} n_{z}}{\theta_{xy}} - \dfrac{2 n_{x} n_{z} n_{xy}}{\theta}\right) (1 - \cos(\theta))$ \\[5mm]
$\dfrac{\partial R_{zy}}{\partial \theta_{xy}}$&$=\left(\dfrac{n_{x}}{\theta_{xy}} - \dfrac{n_{x} n_{xy}}{\theta} + n_{y} n_{z} n_{xy}\right) \sin(\theta) + n_{x} n_{xy} \cos(\theta) + \left(\dfrac{n_{y} n_{z}}{\theta_{xy}} - \dfrac{2 n_{y} n_{z} n_{xy}}{\theta}\right) (1 - \cos(\theta))$ \\[5mm]
$\dfrac{\partial R_{zz}}{\partial \theta_{xy}}$&$=-2n_{xy}n_{z}^2\dfrac{(1-\cos(\theta))}{\theta}-n_{xy}^3(n_{xy}^2+n_{z}^2) \sin(\theta)$ \\[3mm]
\hline
\end{tabular}
\end{table}

\begin{table}[hbt!]
\raggedright
\caption {Derivatives of rotation matrix components with respect to $z$ rotations in polar coordinates of the $xy$-plane for point-to-point optimizations}
\label{tab:pp_awz2} 
\begin{tabular}{| rl |}
  \hline 
& \\
$\dfrac{\partial R_{xx}}{\partial \theta_{z}}$&$=(n_x^2 n_z - n_z)\sin(\theta) - \dfrac{2 n_x^2 n_z}{\theta}(1-\cos(\theta))
$ \\[5mm]
$\dfrac{\partial R_{xy}}{\partial \theta_{z}}$&$=
\Bigg(n_x n_y n_z - \dfrac{\theta_{xy}^2}{\theta^3}\Bigg) \sin(\theta) 
+ \dfrac{2 n_x n_y n_z}{\theta}(\cos(\theta)-1) - n_z^2 \cos(\theta)
$ \\[5mm]
$\dfrac{\partial R_{xz}}{\partial \theta_{z}}$&$=
 \Bigg(n_x n_z^2 - \dfrac{n_y n_z}{\theta} \Bigg) \sin(\theta) + \dfrac{(2 n_x n_z^2 - n_x)}{\theta} (\cos(\theta)-1)  + n_y n_z \cos(\theta)$ \\[5mm]
$\dfrac{\partial R_{yx}}{\partial \theta_{z}}$&$= \Bigg(n_x n_y n_z + \dfrac{\theta_{xy}^2}{\theta^3} \Bigg) \sin(\theta) + \dfrac{2 n_x n_y n_z}{\theta}(\cos(\theta)-1) + n_z^2 \cos(\theta)$ \\[5mm]
$\dfrac{\partial R_{yy}}{\partial \theta_{z}}$&$=(n_y^2 n_z - n_z)\sin(\theta) - \dfrac{2 n_y^2 n_z}{\theta}(1-\cos(\theta))$ \\[5mm]
$\dfrac{\partial R_{yz}}{\partial \theta_{z}}$&$= \Bigg(n_y n_z^2 + \dfrac{n_x n_z}{\theta} \Bigg)\sin(\theta) + \dfrac{(2 n_y n_z^2 - n_y)}{\theta} (\cos(\theta)-1)  - n_x n_z \cos(\theta)$ \\[5mm]
$\dfrac{\partial R_{zx}}{\partial \theta_{z}}$&$= \Bigg(n_x n_z^2 + \dfrac{n_y n_z}{\theta} \Bigg)\sin(\theta) + \dfrac{(2 n_x n_z^2 - n_x)}{\theta} (\cos(\theta)-1)  - n_y n_z \cos(\theta)$ \\[5mm]
$\dfrac{\partial R_{zy}}{\partial \theta_{z}}$&$= \Bigg(n_y n_z^2 - \dfrac{n_x n_z}{\theta} \Bigg)\sin(\theta) + \dfrac{(2 n_y n_z^2 - n_y)}{\theta} (\cos(\theta)-1)  + n_x n_z \cos(\theta)$ \\[5mm]
$\dfrac{\partial R_{zz}}{\partial \theta_{z}}$&$=(n_z^3 - n_z)\sin(\theta) + \dfrac{2(n_z - n_z^3 )}{\theta}(1-\cos(\theta))$ \\[3mm]
\hline
\end{tabular}
\end{table}  

\begin{table}[hbt!]
\raggedright
\caption {Exact gradients of the quaternion elements $A$, $B$, $C$, $D$ with respect to cartesian coordinates ($\theta_x$, $\theta_y$ and $\theta_z$) and polar coordinates ($\alpha$,$\theta_{xy}$ and $\theta_z$). Auxiliary variable is $n_{xy} = \frac{\theta{xy}}{\theta}$, other variables as defined in the main text.
} \label{tab:ur} 
\begin{tabular}{|p{0.7cm}p{7.0 cm}|p{0.7cm}p{8.0 cm}|}
\hline 
cartesian & & polar & \\
  \hline 
& & & \\
$\dfrac{\partial A}{\partial \theta_x}=$ & $ (1-n_{x}^2) \ \dfrac{\sin(\theta/2)}{\theta}+n_{x}^2 \ \dfrac{\cos(\theta/2)}{2}$ & $\dfrac{\partial A}{\partial\alpha}=$ & $- \theta_{xy}\sin(\alpha) \ \dfrac{\sin(\theta/2)}{\theta}$
\\[5mm]
$\dfrac{\partial B}{\partial \theta_x}=$ & $n_{x}n_y\bigg(\dfrac{\cos(\theta/2)}{2}-\dfrac{\sin(\theta/2)}{\theta}\bigg)$ & $\dfrac{\partial B}{\partial\alpha}=$ & $ \ \ \theta_{xy}\cos(\alpha) \ \dfrac{\sin(\theta/2)}{\theta}$
\\[5mm]
$\dfrac{\partial C}{\partial \theta_x}=$ & $n_{x}n_z\bigg(\dfrac{\cos(\theta/2)}{2}-\dfrac{\sin(\theta/2)}{\theta}\bigg)$ & $\dfrac{\partial C}{\partial\alpha}=$ & $0$
\\[5mm]
$\dfrac{\partial D}{\partial \theta_x}=$ & $-n_{x} \ \dfrac{\sin(\theta/2)}{2}$ & $\dfrac{\partial D}{\partial\alpha}=$ & $0$ \\[3mm]
\hline
& & & \\
$\dfrac{\partial A}{\partial \theta_y}=$ & $n_{x}n_y\bigg(\dfrac{\cos(\theta/2)}{2}-\dfrac{\sin(\theta/2)}{\theta}\bigg)$ 
& $\dfrac{\partial A}{\partial\theta_{xy}}=$ & ~~$\cos(\alpha)\bigg( (1-n_{xy}^2) \ \dfrac{\sin(\theta/2)}{\theta}+n_{xy}^2 \ \dfrac{\cos(\theta/2)}{2}\bigg)$ 
\\[5mm]
$\dfrac{\partial B}{\partial \theta_y}=$ & $ (1-n_{y}^2) \ \dfrac{\sin(\theta/2)}{\theta}+n_{y}^2 \ \dfrac{\cos(\theta/2)}{2}$ 
& $\dfrac{\partial B}{\partial\theta_{xy}}=$ & ~~$\sin(\alpha)\bigg( (1-n_{xy}^2) \ \dfrac{\sin(\theta/2)}{\theta}+n_{xy}^2 \ \dfrac{\cos(\theta/2)}{2}\bigg)$
\\[5mm]
$\dfrac{\partial C}{\partial \theta_y}=$ & $n_{y}n_z\bigg(\dfrac{\cos(\theta/2)}{2}-\dfrac{\sin(\theta/2)}{\theta}\bigg)$ 
& $\dfrac{\partial C}{\partial\theta_{xy}}=$ & ~~$n_{xy}n_z\bigg(\dfrac{\cos(\theta/2)}{2}-\dfrac{\sin(\theta/2)}{\theta}\bigg)$ 
\\[5mm]
$\dfrac{\partial D}{\partial \theta_y}=$ & $-n_{y} \ \dfrac{\sin(\theta/2)}{2}$ 
& $\dfrac{\partial D}{\partial\theta_{xy}}=$ & ~~$-n_{xy} \ \dfrac{\sin(\theta/2)}{2}$
\\[3mm]
\hline
& & & \\
$\dfrac{\partial A}{\partial \theta_z}=$ & $n_{x}n_z\bigg(\dfrac{\cos(\theta/2)}{2}-\dfrac{\sin(\theta/2)}{\theta}\bigg)$
& $\dfrac{\partial A}{\partial\theta_{z}}=$ & ~$n_{xy}n_z\cos(\alpha)\bigg(\dfrac{\cos(\theta/2)}{2}-\dfrac{\sin(\theta/2)}{\theta}\bigg)$
\\[5mm]
$\dfrac{\partial B}{\partial \theta_z}=$ & $n_{y}n_z\bigg(\dfrac{\cos(\theta/2)}{2}-\dfrac{\sin(\theta/2)}{\theta}\bigg)$
& $\dfrac{\partial B}{\partial\theta_{z}}=$ & ~$n_{xy}n_z\sin(\alpha)\bigg(\dfrac{\cos(\theta/2)}{2}-\dfrac{\sin(\theta/2)}{\theta}\bigg)$
\\[5mm]
$\dfrac{\partial C}{\partial \theta_z}=$ & $(1-n_z^2) \ \dfrac{\sin(\theta/2)}{\theta}+n_z^2 \ \dfrac{\cos(\theta/2)}{2}$
& $\dfrac{\partial C}{\partial\theta_{z}}=$ & $(1-n_z^2) \ \dfrac{\sin(\theta/2)}{\theta}+n_z^2 \ \dfrac{\cos(\theta/2)}{2}$
\\[5mm]
$\dfrac{\partial D}{\partial \theta_z}=$ & $-n_z \ \dfrac{\sin(\theta/2)}{2}$ & 
$\dfrac{\partial D}{\partial\theta_{z}}=$ & $-n_z \ \dfrac{\sin(\theta/2)}{2}$
\\[3mm]
\hline
\end{tabular}
\end{table}

\begin{table}[hbt!]
\caption {Benchmark runtimes of 1000 exact gradient calculations using the different formulae from Tables \ref{tab:pp_dx}-\ref{tab:pp_awz1}} \label{tab:runt} 
\begin{tabular}{crrcrr}
\hline
& & & & & \\
PP (Cartesian) &  Laptop$^a$         & workstation$^b$      & UR (Cartesian) &    Laptop$^a$        & workstation$^b$ \\[3mm]
\hline
\hline
analytical              & & & analytical     & & \\[1mm]
$\partial R / \partial \theta_x$ &   325.9 $\mu$s &   368.0 $\mu$s   & $\partial Q / \partial \theta_x$ &     271.4 $\mu$s  &    331.3 $\mu$s  \\
$\partial R / \partial \theta_y$ &   298.2 $\mu$s &   361.7 $\mu$s   & $\partial Q / \partial \theta_y$ &     278.4 $\mu$s  &    329.7 $\mu$s  \\
$\partial R / \partial \theta_z$ &   326.0 $\mu$s &   362.6 $\mu$s   & $\partial Q / \partial \theta_z$ &     288.7 $\mu$s  &    333.4 $\mu$s  \\[3mm]
exponential               & & & exponential     & & \\[1mm]
$\partial R / \partial \theta_x$ & 45788.0 $\mu$s &   33467.0 $\mu$s & $\partial Q / \partial \theta_x$ &    43774.0 $\mu$s &    32058.0 $\mu$s  \\ 
$\partial R / \partial \theta_y$ & 49442.0 $\mu$s &   33990.0 $\mu$s & $\partial Q / \partial \theta_y$ &    49276.0 $\mu$s &    32492.0 $\mu$s  \\ 
$\partial R / \partial \theta_z$ & 51196.0 $\mu$s &   33152.0 $\mu$s & $\partial Q / \partial \theta_z$ &    48715.0 $\mu$s &    32651.0 $\mu$s  \\ [3mm]
finite differences               & & & finite differences     & & \\[1mm]
$\partial R / \partial \theta_x$ &   490.8 $\mu$s &   915.5 $\mu$s   & $\partial Q / \partial \theta_x$ &     311.3 $\mu$s  &    328.9 $\mu$s  \\
$\partial R / \partial \theta_y$ &   505.4 $\mu$s &   907.5 $\mu$s   & $\partial Q / \partial \theta_y$ &     302.7 $\mu$s  &    325.3 $\mu$s  \\
$\partial R / \partial \theta_z$ &   501.5 $\mu$s &   908.0 $\mu$s   & $\partial Q / \partial \theta_z$ &     307.3 $\mu$s  &    325.3 $\mu$s  \\[3mm]
\hline
& & & & & \\
PP (polar)   &    &     &    UR (polar)   &    &   \\[3mm]
\hline
\hline
analytical              & & & analytical     & & \\[1mm]
$\partial R / \partial \alpha$  &   286.6 $\mu$s &   337.4 $\mu$s   & $\partial R / \partial \alpha$ &     287.4 $\mu$s   &   330.2 $\mu$s    \\
$\partial R / \partial \theta_{xy}$ &   293.5 $\mu$s &   349.2 $\mu$s   & $\partial R / \partial \theta_{xy}$ &     281.9 $\mu$s   &   354.0 $\mu$s    \\
$\partial R / \partial \theta_z$ &   276.1 $\mu$s &   335.5 $\mu$s   &$\partial R / \partial \theta_{z}$ &      277.0 $\mu$s   &   332.5 $\mu$s    \\[3mm]
finite differences             & & & finite differences     & & \\[1mm]
$\partial R / \partial \alpha$ &   453.1 $\mu$s &   435.8 $\mu$s   & $\partial R / \partial \alpha$ &     342.3 $\mu$s   &   407.3 $\mu$s    \\
$\partial R / \partial \theta_{xy}$ &   451.6 $\mu$s &   433.9 $\mu$s   & $\partial R / \partial \theta_{xy}$ &     410.9 $\mu$s   &   412.9 $\mu$s   \\
$\partial R / \partial \theta_z$ &   457.3 $\mu$s &   443.6 $\mu$s   & $\partial R / \partial \theta_{z}$ &     412.7 $\mu$s   &   416.0 $\mu$s   \\[3mm]
\hline
\end{tabular}

$^a$  Lenovo Thinkpad X1 (2023) 12th Gen Intel Core i7-1260P 2.10 GHz; $^b$ AMD Ryzen 9 5900X 12-Core Processor 3.70 GHz.
\end{table}

\begin{table}[hbt!]
\caption {Summary of individual optimizations performed using the different formulae from Tables \ref{tab:pp_dx}-\ref{tab:pp_awz1}. In all cases, pulses of 500~$\mu$s duration, maximum rf-amplitude of 5~kHz with $\pm$10\% B$_1$-compensation (3 $B_1$ points), and a bandwidth of 6000~kHz (11 offset points, covering 50~ppm in $^{15}$N at a 1.2~GHz spectrometer) were optimized using a BFGS optimization algorithm on a single core of a laptop.  } \label{tab:15N} 
\begin{tabular}{lccccc}
\hline
& & & & & \\
Controls &  $\Delta$ t &    iterations  & optimization & time per & quality \\
 &  [$\mu$s]  & \# & time [s] & iter. [$\mu$s] &  factor\\[3mm]
\hline
\hline
excitation  & & & & & \\[1mm]
 $\theta_x$, $\theta_y$					   & 1 / 10 / 50 & 	4944 / 424  / 148 & 665,4 / 7,4  / 0,598 & 134,5 / 17,4 / 4,0 & 0,9987 / 0,9960 / 0,9973 \\
 $\theta_x$, $\theta_y$, $\theta_z$		   & 1 / 10 / 50 & 	1294 / 152  / 98  & 306,2 / 7,2  / 0,277 & 236,6 / 47,3 / 2,8 & 0,9990 / 0,9938 / 0,9955 \\
 $\theta_{xy}$, $\alpha$				   & 1 / 10 / 50 & 	777  / 134  / 119 & 383,9 / 7,6  / 0,608 & 494,0 / 56,7 / 5,1 & 0,9962 / 0,9949 / {\it 0,9923} \\
 $\theta_{xy}$, $\alpha$, $\theta_z$	   & 1 / 10 / 50 & 	1083 / 198  / 65  & 991,3 / 12,8 / 0,428 & 915,3 / 64,6 / 6,5 & 0,9978 / 0,9981 / 0,9968 \\
 $\theta_{xy}^{\rm{red}}$, $\alpha$			   & 1 / 10 / 50 & 	147  / 146  / 105 & 17,7  / 2,8  / 0,407 & 120,4 / 19,1 / 3,8 & 0,9966 / 0,9991 / 0,9985 \\
 $\theta_{xy}^{\rm{red}}$, $\alpha$, $\theta_z$ & 1 / 10 / 50 & 	765  / 247  / 85  & 82,7  / 5,5  / 0,421 & 108,1 / 22,2 / 4,9 & 0,9989 / 0,9989 / 0,9968 \\
 $\alpha$								   & 1 / 10 / 50 & 	847  / 1018 / 82  & 53,7  / 8,1  / 0,126 & 63,4  / 7,9  / 1,5 & 0,9981 / {\it 0,9995} / 0,9963 \\[3mm]
inversion   & & & & & \\[1mm]
 $\theta_x$, $\theta_y$				&  1 / 10 / 50  &  4805 / 231  / 87  & 727.4 / 5.4  / 0.277  & 151.3 / 23.3 / 3.1 & 0.9978 / 0.9947 / 0.9926 \\
 $\theta_x$, $\theta_y$, $\theta_z$		&  1 / 10 / 50  &  1427 / 185  / 70  & 316.7 / 6.8  / 0.418  & 221.9 / 36.7 / 5.9 & 0.9972 / 0.9971 / 0.9946 \\
 $\theta_{xy}$, $\alpha$				&  1 / 10 / 50  &  809  / 1192 / 78  & 488.6 / 21   / 0.318  & 603.9 / 17.6 / 4.0 & {\it 0.9926} / 0.9996 / 0.9972 \\
 $\theta_{xy}$, $\alpha$, $\theta_z$		&  1 / 10 / 50  &  819  / 173  / 132 & 632.9 / 11.3 / 0.74   & 772.7 / 65.3 / 5.6 & 0.9950 / 0.9961 / 0.9987 \\
 $\theta_{xy}^{\rm{red}}$, $\alpha$			&  1 / 10 / 50  &  331  / 302  / 95  & 46.2  / 3.5  / 0.376  & 139.5 / 11.5 / 3.9 & 0.9972 / 0.9995 / 0.9973 \\
 $\theta_{xy}^{\rm{red}}$, $\alpha$, $\theta_z$	&  1 / 10 / 50  &  260  / 121  / 54  & 40.8  / 2.5  / 0.228  & 156.9 / 20.6 / 4.2 & 0.9965 / 0.9995 / 0.9926 \\
 $\alpha$						&  1 / 10 / 50  &  573  / 753  / 88  & 39.2  / 4.4  / 0.215  & 68.4  / 5.8  / 2.4 & {\it 0.9997} / 0.9981 / 0.9932 \\[3mm]
\hline
\end{tabular}
\end{table}

\begin{table}[hbt!]
\caption {Different pulse types optimized for a typical $^{13}$C scenario, where 250~ppm on a 600~MHz spectrometer or 160~ppm on a 1.0 GHz spectrometer have to be covered. Pulses were optimized using an L-BFGS optimization algorithm on a single core of a laptop.  } \label{tab:13C} 
\begin{tabular}{lccccccccccc}
\hline
& & & & & & & & & & & \\
Pulse & Controls & offset & BW   & B$_1$ & $\pm \vartheta$ & $\Delta$ t & $t_p$ & $\sqrt{\overline{P}^{\rm{max}}}$ / & iter.  & opt.       & quality \\
type  &  &  \#    & [kHz] & \#     &       [\%]          &[$\mu$s]  & [$\mu$s]  & $\theta_{xy}^{\rm{max}}$ [kHz]     & \#    & time [s] &  factor\\[3mm]
\hline
\hline
excitation & $\theta_x$, $\theta_y$			& 31 & 40 & 3 & 5 & 1  & 500  & 10/20 & 3956 & 1664.0   & {\it 0.99982} \\
excitation & $\theta_x$, $\theta_y$, $\theta_z$	& 31 & 40 & 3 & 5 & 10 & 500  & 10/20 & 90   & 6.4    & 0.99914 \\
excitation & $\alpha$					& 31 & 40 & 3 & 5 & 1  & 500  & - /10 & 1248 & 293.0  & {\it 0.99715} \\
inversion  & $\theta_x$, $\theta_y$			& 31 & 40 & 3 & 5 & 1  & 500  & 10/20 & 1966 & 690.0  & {\it 0.99983} \\
inversion  & $\theta_x$, $\theta_y$, $\theta_z$	& 31 & 40 & 3 & 5 & 10 & 500  & 10/20 & 104  & 7.5    & 0.99839 \\
inversion  & $\alpha$					& 31 & 40 & 3 & 5 & 1  & 500  & - /10 & 1473 & 344.0  & {\it 0.99744}  \\ 
UR-90$^{\circ}$  & $\theta_x$, $\theta_y$		& 61 & 40 & 3 & 5 & 1  & 1000 & 10/20 & 5000$^*$ & 3170.0 & {\it 0.99996} \\ 
UR-180$^{\circ}$ & $\theta_x$, $\theta_y$		& 61 & 40 & 3 & 5 & 1  & 1000 & 10/20 & 2854 & 3353.6 & {\it 0.99998} \\[3mm]

\hline
\end{tabular}

$^*$ optimization stopped after reaching maximum number of iterations.
\end{table}

\begin{table}[hbt!]
\caption {Different pulse types optimized for a typical $^{19}$F scenario, where 200~ppm need to be covered on a 600~MHz spectrometer. Pulses were optimized using an L-BFGS optimization algorithm on a single core of a laptop.  } \label{tab:19F} 
\begin{tabular}{lccccccccccc}
\hline
& & & & & & & & & & & \\
Pulse & Controls & offset & BW   & B$_1$ & $\pm \vartheta$ & $\Delta$ t & $t_p$ & $\sqrt{\overline{P}^{\rm{max}}}$ / & iter.  & opt.       & quality \\
type  &  &  \#    & [kHz] & \#     &       [\%]          &[$\mu$s]  & [$\mu$s]  & $\theta_{xy}^{\rm{max}}$ [kHz]     & \#    & time [s] &  factor\\[3mm]
\hline
\hline
saturation & $\alpha$					& 31  & 120 & 1 & -    & 2  & 120  & -/10 & 244  & 1.6   & {\it 0.99999} \\
excitation & $\theta_{xy}^{\rm{red}}$, $\alpha$	& 301 & 120 & 5 & 10 & 2 & 1400  & 10/- & 5000$^*$   & 94570    & {\it 0.99627} \\
excitation & $\alpha$					& 121 & 120 & 5 & 10 & 1  & 350   & -/20 & 894 & 904.4  & {\it 0.99419} \\
inversion  & $\theta_{xy}^{\rm{red}}$, $\alpha$	& 301 & 120 & 5 & 10 & 2 & 1400  & 10/- & 5000$^*$  & 42864    & {\it 0.99989} \\
inversion  & $\alpha$					& 121 & 120 & 5 & 10 & 1  & 300  & -/20 & 567 & 757.1  & {\it 0.99885}  \\ 
UR-90$^{\circ}$  & $\alpha$				& 121 & 120 & 5 & 10 & 1  & 600 & -/20 & 2587 & 2833.2 & {\it 0.99915} \\ 
UR-180$^{\circ}$ & $\alpha$				& 121 & 120 & 5 & 10 & 1  & 700 & -/20 & 10000$^*$  & 15411.0 & {\it 0.99951} \\[3mm]

\hline
\end{tabular}

$^*$ optimization stopped after reaching maximum number of iterations.
\end{table}

\end{document}